\documentstyle[12pt]{article}
\begin{document}
\author{S.V. Ludkovsky.}
\title{Approximate and global differentiability of functions
over non-archimedean fields.}
\date{23.08.2006}
\maketitle

\begin{abstract}
The article is devoted to the investigation of approximate, global
and along curves smoothness of functions $f(x_1,...,x_m)$ of
variables $x_1,...,x_m$ in infinite fields with non trivial
multiplicative ultra-norms and relations between them.  Then classes
of smoothness $C^{n,r}$ and $C^{n,r}_b$ and more general in the
sense of Lipschitz for partial difference quotients are considered
and theorems for them are proved. Moreover, an approximate
differentiability of functions relative to measures is defined and
investigated. Its relations with lipschitzian property and almost
everywhere differentiability are studied. Finally theorems about
relations between approximate differentiability by all variables and
along curves are proved.
\end{abstract}
\section{Introduction}
\par Fields with non-archimedean valuations such as the field of
$p$-adic numbers were first introduced by K. Hensel \cite{hensel}.
Then it was proved by A. Ostrowski \cite{ostrow} that on the field
of rational numbers each multiplicative norm is either the usual
norm as in $\bf R$ or is equivalent to a non-archimedean norm
$|x|=p^{-k}$, where $x=np^k/m\in \bf Q$, $n, m, k\in \bf Z$, $p\ge
2$ is a prime number, $n$ and $m$ and $p$ are mutually pairwise
prime numbers. It is well known, that each locally compact infinite
field with a non trivial non-archimedean valuation is either a
finite algebraic extension of the field of $p$-adic numbers or is
isomorphic to the field ${\bf F}_{p^k}(\theta )$ of power series of
the variable $\theta $ with expansion coefficients in the finite
field ${\bf F}_{p^k}$ of $p^k$ elements, where $p\ge 2$ is a prime
number, $k\in \bf N$ is a natural number \cite{roo,weil}. Non
locally compact fields are also wide spread \cite{diar,roo,sch1}.
\par Last years non-archimedean analysis \cite{roo,sch1,sch2}
and mathematical physics \cite{yojan,khrhs,khrqp,vla3,padicmp} are
being fastly delevoped. But many questions and problems remain open.
\par In the non-archimedean analysis classes of smoothness are defined
in another fashion as in the classical case over $\bf R$, since
locally constant functions on fields $\bf K$ with non-archimedean
valuations are infinite differentiable and there exist non trivial
non locally constant functions infinite differentiable with
identically zero derivatives \cite{sch1,sch2}. This is caused by the
stronger ultrametric inequality $|x+y|\le \max (|x|,|y|)$ in
comparison with the usual triangle inequality, where $|x|$ is a
multiplicative norm in $\bf K$ \cite{roo}. In papers
\cite{luseamb,lutmf99,luanmbp,beglneeb} there were considered
classes of smoothness $C^n$ for functions of several variables in
non-archimedean fields or in topological vector spaces over such
fields.
\par In the classical functional analysis the approximate
differentiability and almost everywhere differentiability are widely
used and studied, for example, for the needs of the geometric
measure theory \cite{federer}. On the other hand, in the
non-archimedean case this subject was not yet investigated.
\par This paper is devoted to
the investigation of smoothness of functions $f(x_1,...,x_m)$ of
variables $x_1,...,x_m$ in infinite fields with non trivial
non-archimedean valuations, where $m\ge 2$. In the paper fields
locally compact and as well as non locally compact are considered.
Theorems about classes of smoothness $C^n$ or $C^n_b$ of functions
with continuous or bounded uniformly continuous on bounded domains
partial difference quotients up to the order $n$ are investigated.
\par In the second section the approximate limits and approximate
differentiability in the sense of partial difference quotients are
defined and investigated over locally compact fields relative to the
Haar nonnegative measures on fields. Non-archimedean analogs of
classical theorems of Kirzsbraun, Rademacher, Stepanoff, Whitney are
formulated and proved (see Theorems 2.8, 9, 17, 22 respectively).
Their relations with the lipschitzian property and almost everywhere
differentiability are studied (see Theorems 2.15, 19, 20, 23, 24).
Finally theorems 2.25 and 26 about relations between approximate
differentiability by all variables and along curves are proved.
Frequently formulations of theorems and their proofs in the
non-archimedean case are different from the classical case due to
specific features of the non-archimedean analysis and underlying
fields. All results of this paper are obtained for the first time.

\section{Approximate differentiability of functions}
\par {\bf 1. Definition and Notations.} Let $\bf K$ be an infinite
locally compact field with a non trivial non archimedean valuation
and ${\cal B}({\bf K})={\cal B}$ be a Borel $\sigma $-algebra of
subsets of $\bf K$. A $\sigma $-additive $\sigma $-finite measure
$\mu : {\cal B}\to [0,\infty ]$ is called the Haar measure, if $\mu
$ is non zero and $\mu (x+A)=\mu (A)$ for each $x\in \bf K$ and
$A\in {\cal B}$. For convenience put $\mu (B({\bf K},0,1))=1$ and
choose an equivalent valuation $|x|=|x|_{\bf K}=mod_{\bf K}(x)$ in
$\bf K$, where $mod_{\bf K}(x)$ is the modular function such that
$\mu (x B({\bf K},0,R))=mod_{\bf K}(x)\mu (B({\bf K},0,R))$ for each
$x\in \bf K$, where $R$ belongs to the valuation group $\Gamma _{\bf
K}:= \{ |x|: 0\ne x\in {\bf K} \} $ of $\bf K$, $A_1A_2:= \{ a:
a=a_1a_2, a_1\in A_1, a_2\in A_2 \} $ for $A_1, A_2\subset \bf K$,
$B({\bf K}^m,y,R) := \{ z\in {\bf K}^m: |z-y|\le R \} $, $m\in \bf
N$. For ${\bf K}^m$ take the measure $\mu ^m=\otimes_{j=1}^m\mu _j$,
where $\mu _j=\mu =:\mu ^1$ for each $j$, such that $\mu
^m(A_1\times ... \times A_m)=\prod_{j=1}^m \mu (A_j)$ for each
$A_1,...,A_m\in {\cal B}({\bf K})$, $\otimes_{j=1}^m{\cal B}({\bf
K})$ is the minimal $\sigma $-algebra generated by subsets of the
form $A_1\times ...\times A_m$ and the Borel $\sigma $-algebra
${\cal B}({\bf K}^m)$ of ${\bf K}^m$ coincides with
$\otimes_{j=1}^m{\cal B}({\bf K})$.
\par {\bf 2. Remark.} Suppose
that $(X,\rho _X)$ is a metric space with a set $X$ and a metric
$\rho _X$ in it. A non negative measure $\nu $ on a $\sigma
$-algebra $\cal C$ of $X$ is called Borel regular if and only if
${\cal B}(X)\subset \cal C$ and for each $A\in \cal C$ there exists
$H\in {\cal B}(X)$ such that $\nu (H)=\nu (A)$, where ${\cal B}(X)$
denotes the Borel $\sigma $-algebra of $X$ which is the minimal
$\sigma $-algebra generated by open subsets of $X$. Denote by $\cal
M$ the class of all Borel regular non negative measures $\nu $ on
$X$ such that each bounded subset $A$ in $X$ has a finite measure
$0\le \nu (A)<\infty $.
\par Consider the family $\cal N$ of all subsets $A$ of $X$
for which there exists $G=G(A)\in {\cal B}(X)$ such that $A\subset
G$ and $\nu (G)=0$. The minimal $\sigma $-algebra ${\cal A}_{\nu
}={\cal A}_{\nu }(X)$ generated by ${\cal B}(X)\cup \cal N$ is the
$\nu $-completion of ${\cal B}(X)$ and it consists of $\nu
$-measurable subsets.
\par A subset of the form $ \{ (y,A): y\in A\subset X \} $ is called
a covering relation. If $Y$ is a subset of $X$, then put $V(Y):= \{
A: \mbox{ there exists } y\in Y, (y,A)\in V \} $, where $V$ is a
covering relation in $X$. Then $V$ is called fine at a point $y$ if
and only if $\inf_{(y,A)\in V} diam (A)=0$, where with each
$A\subset X$ is associated its diameter $diam (A):=\sup_{x, y\in A}
\rho _X(x,y)$.
\par Consider a covering relation $V$ in $X$ satisfying conditions:
\par $(V1)$ $V(X)$ is a family of Borel subsets of $X$,
\par $(V2)$ $V$ is fine at each point of $X$,
\par $(V3)$ if $C\subset V$ and $Y\subset X$ and $C$ is fine at each
point of $Y$, then $C(Y)$ has a countable disjoint subfamily
covering $\nu $ almost all of $Y$.
\par If for a given $\nu \in \cal M$ a covering relation
$V$ satisfies Conditions $(V1-V3)$, then it is called a $\nu $
Vitali relation.
\par Suppose that $\nu \in \cal M$ and and $V$ is a $\nu $ Vitali
relation. With each $\lambda \in \cal M$ there is associated another
measure $\lambda _{\nu }$ defined by the formula:
\par $\lambda _{\nu }(A) := \inf \{ \lambda (S): S \in {\cal B}(X),
\nu (A\setminus S)=0 \} $, whenever $A\subset X$, hence $\lambda
_{\nu }\le \lambda $.
\par If a covering relation $V$ is fine at a point $x\in X$
and $f: dom (f) \to [-\infty , \infty ]$, where $dom (f)\subset
{\cal B}(X)$ is the domain of $f$, then
\par $(V) \lim_{S\to x} f(S) := \lim_{0<\epsilon \to 0}
\{ f(S): (x,S)\in V, diam (S)<\epsilon , S\in dom (f) \} $, \\
analogously are defined \par $(V) \lim_{S\to x}\sup
f(S):=\lim_{0<\epsilon \to 0} \sup \{ f(S): (x,S)\in V, diam
(S)<\epsilon , S\in dom (f) \} $ \\ and $(V) \lim_{S\to x} \inf f$.
\par For a subset $A$ in $X$ and a point $x\in X$ the limit
$(V) \lim_{S\to x} \nu (S\cap A)/\nu (S)$ is called the $(\nu ,V)$
density of $A$ at $x$.
\par If $g: X\to Y$ is a mapping of a metric space $(X, \rho _X)$
into a (Hausdorff) topological space $Y$, then $y\in Y$ is called an
approximate limit of $g$ at $x$ (relative to a measure $\nu \in
{\cal M}(X)$ and a $\nu $ Vitaly relation $V$) if and only of for
each neighborhood $W$ of $y$ in $Y$ the set $X\setminus g^{-1}(W)$
has zero density at $x$ and it is denoted by $y=(\nu ,V) ap
\lim_{z\to x} g(z)$. If $(\nu ,V)$ are specified, then they may be
omitted for brevity.
\par A function $g$ is called $(\nu ,V)$ approximately continuous if
and only if $x\in dom (g)$ and $(\nu ,V) ap \lim_{z\to x}
g(z)=g(x)$.
\par {\bf 3. Definitions.} Let $Y$ be a topological vector space
over $\bf K$, $g: U\to Y$ be a mapping, where $U$ is open in ${\bf
K}^m$, $m\in \bf N$. Then $g$ is called approximate differentiable
at a point $x$ of $U$ if there exists an open neighborhood $W$ of
$x$, $W\subset U$, such that ${\bar {\Phi }}^1 g(x;v;t)$ is $\mu
^{2m+1}$ almost everywhere continuous on $W^{(1)}$ and at $x$ there
exists a linear mapping $T: {\bf K}^m\to Y$ such that $(\mu ^m,V) ap
\lim_{z\to x} |g(z)-g(x)- T (z-x)|/ |z-x|=0$, where $V = {\cal
B}({\bf K}^m)$. This $T$ is also denoted by $\mbox{ }_{ap} Dg(x)$.
If this is satisfied for each $z\in U$, then $g$ is called
approximate differentiable on $U$. The family of all such functions
denote by $\mbox{ }_{ap}C^1(U,Y)$. Then also define approximate
partial derivatives:
\par $\mbox{ }_{ap} D_jg(x) := ap \lim_{t\to 0} [g(x_1,...,
x_{j-1},x_j+t,x_{j+1},...,x_m)-g(x)]/t$.
\par The family of all $f\in C^n(U,Y)$ such that ${\bar {\Phi }}^nf
\in \mbox{ }_{ap}C^1(U^{(n)},Y)$ denote by $\mbox{
}_{ap}C^{n+1}(U,Y)$. Suppose now that $A$ is a $\mu ^m$ measurable
subset of $U$ and $Y$ be a normed space, $0<r\le 1$, then denote by
$\mbox{ }_{ap}C^{n,r}(U,A,Y)$ the family of all $f\in C^n(U,Y)$ with
\par $ap {\overline {\lim}}_{x^{(n)}\to z^{(n)}} \| {\bar {\Phi
}}^kf(x^{(n)})-{\bar {\Phi }}^nf(z^{(n)}) \|
_{C^0(V^{(n)}_{z,R},Y)}/|x^{(n)}-z^{(n)}|^r<\infty $ \\
for $\mu ^m$ almost all $z\in A$ and each $0<R<\infty $, where
$V^{(k)}_{z,R}$ corresponds to $U_{z,R}=U\cap B({\bf K}^m,z,R)$.
\par {\bf 4. Lemma.} {\it The families $\mbox{ }_{ap}C^{n+1}(U,Y)$
and $\mbox{ }_{ap}C^{n,r}(U,A,Y)$ are the $\bf K$-linear spaces.}
\par {\bf Proof.} 1. Since $C^n(U,Y)$ is the $\bf K$ linear space, then
it is sufficient to verify, that the set of all ${\bar {\Phi }}^nf$
with $f\in \mbox{ }_{ap}C^{n+1}(U,Y)$ is $\bf K$ linear. Therefore,
the consideration reduces to $\mbox{ }_{ap}C^1(U,Y)$, where
$U^{(n)}$ is denoted also by $U$. If $f, g\in \mbox{ }_{ap}C^1(U,Y)$
and $a, b\in \bf K$, then
\par ${\bar {\Phi }}^1(af+bg)(x;v;t)=a{\bar {\Phi }}^1f(x;v;t)+
b{\bar {\Phi }}^1g(x;v;t)$,\\ since for $\mu ^{2m+1}$-almost all
points $(x;v;t)\in {\bf K}^{2m+1}$ the right side terms are
continuous and hence the left side term is such also. If $x\in U$,
then
\par $ap \lim_{z\to x}|(af+bg)(z)-(af+bg)(x) - a\mbox{
}_{ap}Df(x).(z-x)-b\mbox{ }_{ap}Dg(x).(z-x)|/ |z-x|\le ap \lim_{z\to
x}|a| |f(z)-f(x) - \mbox{ }_{ap}Df(x).(z-x)|/ |z-x| + ap \lim_{z\to
x} |b| |g(z)-g(x) - \mbox{ }_{ap}Dg(x).(z-x)|/ |z-x| =0$, hence
there exists $\mbox{ }_{ap}D(af+bg)(x)=a\mbox{ }_{ap}Df(x)+b\mbox{
}_{ap}Dg(x)$. \par For each $f\in \mbox{ }_{ap}C^1(U,Y)$ the
operator $T=\mbox{ }_{ap}Df(x)$ is unique, since the difference
$H=T_1-T_2$ of two such $\bf K$-linear mappings is subordinated to
the condition: $ap \lim_{v\to 0} |Hv|/|v|=ap \lim_{z\to a}
|H(z-x)|/|z-x|=0$ due to Definition 3. Therefore, if $0<\epsilon
<1$, then there exists $R>0$, $R\in \Gamma _{\bf K}$, such that $\mu
^m(B({\bf K}^m,0,R)\cap \{ v\in {\bf K}^m: |Hv|>\epsilon |v| \} )<
\epsilon ^mR^m$. If $w\in B({\bf K}^m,0,R)$ and $v\in B({\bf
K}^m,w,\epsilon R)$ with $|Hv|\le \epsilon |v|$, then $|Hw|\le \max
(|H(w-v)|, |Hv|)\le \epsilon R\max (\| H \|, 1)$, consequently, $\|
H \| \le \epsilon \max (1, \| H \| )$, consequently, $\| H \| =0$,
since $\epsilon >0$ can be chosen arbitrary small. At the same time
${\bar {\Phi }}^1f(x;v;t)$ is unique on $U^{(1)}$ up to a set of
$\mu ^{2m+1}$-measure zero.
\par 2. The second assertion follows from the inequality
\par $ap {\overline {\lim}}_{x^{(n)}\to z^{(n)}} \| {\bar {\Phi
}}^k(af+bg)(x^{(n)})-{\bar {\Phi }}^n(af+bg)(z^{(n)}) \|
_{C^0(V^{(n)}_{z,R},Y)}/|x^{(n)}-z^{(n)}|^r \le \max (ap {\overline
{\lim}}_{x^{(n)}\to z^{(n)}} |a| \| {\bar {\Phi }}^kf(x^{(n)})-{\bar
{\Phi }}^nf(z^{(n)}) \|
_{C^0(V^{(n)}_{z,R},Y)}/|x^{(n)}-z^{(n)}|^r;$ \\  $|b| ap {\overline
{\lim}}_{x^{(n)}\to z^{(n)}} \| {\bar {\Phi }}^kg(x^{(n)})-{\bar
{\Phi }}^ng(z^{(n)}) \|
_{C^0(V^{(n)}_{z,R},Y)}/|x^{(n)}-z^{(n)}|^r<\infty $ \\
for each $a, b\in \bf K$ and $f, g \in \mbox{ }_{ap}C^{n,r}(U,A,Y).$
\par {\bf 5. Note.} For a locally compact field $\bf K$ each
$x\in \bf K$ has the decomposition $x=\sum_na_n\pi ^n$ for $char
({\bf K})=0$ and $x=\sum_na_n\theta ^n$ for $char ({\bf K})=p>0$
(see Introduction), where $a_n=a_n(x)$ are expansion coefficients.
Introduce on $\bf K$ the linear ordering: $x \prec y$ if and only if
there exists $m\in \bf Z$ such that $a_n(x)=a_n(y)$ for each $n<m$
and $a_m(x)<a_m(y)$ with the natural ordering either in $B({\bf
K},0,1)/B({\bf K},0,|\pi |)$ or in ${\bf F}_{p^k}$ respectively. In
the contrary case we write $x=y$. This linear ordering is compatible
with neither the additive nor the multiplicative structure of $\bf
K$, but it is useful and it was introduced by M. van der Put
\cite{sch1}. In ${\bf K}^m$ we can consider the linear ordering:
$x\prec y$ if and only if there exists $l\in \bf N$ such that
$\mbox{ }_jx=\mbox{ }_jy$ for each $1\le j<l$ and $\mbox{ }_lx\prec
\mbox{ }_ly$, where $x=(\mbox{ }_1x,...,\mbox{ }_mx)\in {\bf K}^m$,
$\mbox{ }_ix\in \bf K$ for each $i=1,...,m$.
\par {\bf 6. Theorem.} {\it Let $\nu $ be a measure on ${\bf K}^m$,
$\nu \in {\cal M}({\bf K}^m)$, let also $f: {\bf K}^m\to \bf K$ be
$\nu $-measurable and $y_1, y_2, y_3,...\in {\bf K}$ are pairwise
distinct and such that for each $y\in {\bf K}$ and every $\epsilon
>0$ there exist $n, k\in \bf N$ such that $|y_n-y|< \epsilon $ and
$y_n\preceq y$ and $|y_k-y|<\epsilon $ and $y\preceq y_k$. Then
there exist $\nu $-measurable subsets $A_1, A_2, A_3,...$ in ${\bf
K}^m$ with characteristic functions $g_j=ch_{A_j}$ such that
\par $f(x)=\sum_{n=1}^{\infty }y_ng_n(x)$ \\
for each $x\in {\bf K}^m$.}
\par {\bf Proof.} If $|y_l-y_j|\le \delta $,
where $0<\delta $, then $\max (|y-y_l|, |y-y_j|)\le \delta $ for
each $y_l\preceq y\preceq y_j$. In the space $L^{\infty }({\bf
K}^m,\nu ,{\bf K})$ the $\bf K$-linear span of characteristic
functions of clopen subsets is dense. Consider the set $W_n:= \{
x\in {\bf K}: |f(x)|\le n \} $, where $n\in \bf N$, then
$f|_{W_n}\in L^{\infty }({\bf K}^m,\nu ,{\bf K})$. Therefore,
$f(x)|_{W_n\setminus W_{n-1}}=\sum_{k=1}^{\infty
}s_{n,k}g_{n,k}(x)$, where each $g_{n,k}$ is the characteristic
function of a $\nu $-measurable subset in ${\bf K}^m$, $s_{n,k}\in
\bf K$, $W_0:=\emptyset $. Then $f(x)=\sum_{k,n=1}^{\infty }
s_{n,k}g_{n,k}(x)$ converges pointwise for each $x\in {\bf K}^m$.
\par Therefore, it is enough to construct the decomposition of
$f|_{W_q\setminus W_{q-1}}$ for arbitrary $q$. Thus consider the
subset $y_n$ with $|y_n|\le q$. For each $\epsilon >0$ there exists
a finite $\epsilon $-net $\{ y_{l(s)}: s=1,..., b \} $,
$b=b(\epsilon )$, that is for each $y\in {\bf K}$ there exists $s$,
$1\le s\le b$, such that $|y-y_{l(s)}|<\epsilon $.
\par Take the sequence $\epsilon _j=|\pi |^j$ and $b_j=b(\epsilon
_j)$. If $y_n$ is the last point of the $|\pi |^j$ net, then new
$|\pi |^{j+1}$-net begins and put $u=1$, if it is not so, then take
$u=0$. Suppose that $\| f-\sum_{j=1}^{n-1}y_jg_j \| _{L^{\infty
}(A_n,\nu ,{\bf K})}>0$. Otherwise the decomposition is already
found. Put by induction $A_n=A_{n,n+l} := \{ x: y_{n}\preceq
f(x)-\sum_{j=1}^{n-1} y_jg_j(x)\prec y_{n+l} \} $ if $y_{n}\preceq
y_{n-k}$, where $l\ge 1$ is the minimal natural number for which
$y_{n}\prec y_{n+l}$ and $|\pi |^{u+1} \| f(x)-\sum_{j=1}^{n-1}
y_jg_j(x)\|_{L^{\infty }(A_n,\nu ,{\bf K})}<|y_{n}-y_{n+l}|\le |\pi
|^u \| f(x)-\sum_{j=1}^{n-1} y_jg_j(x)\|_{L^{\infty }(A_n,\nu ,{\bf
K})}$; $A_{n}=A_{n,n+l} := \{ x: y_{n+l}\preceq
f(x)-\sum_{j=1}^{n-1} y_jg_j(x)\prec y_{n} \} $ if $y_{n-k}\preceq
y_{n}$, where $l\ge 1$ is the minimal natural number for which
$y_{n+l}\prec y_{n}$ and $|\pi |^{u+1} \| f(x)-\sum_{j=1}^{n-1}
y_jg_j(x)\|_{L^{\infty }(A_n,\nu ,{\bf K})}< |y_{n}-y_{n+l}|\le |\pi
|^u \| f(x)-\sum_{j=1}^{n-1} y_jg_j(x)\|_{L^{\infty }(A_n,\nu ,{\bf
K})}$, where $k=k(n)$ is a gap on the preceding step as $l$ on this
step. In accordance with this algorithm some $A_q$ may be empty,
when $n<q<n+l$ for subsequent numbers of the algorithm, so that
$g_q=0$ for such $q$. Then consider the sum
$(\sum_{j=1}^{n-1}y_jg_j(x))+y_{n}g_{n}$. Since $|y_n-y_{n+l}|\le
|\pi |^u |y_n-y_{n-k}|$ for each $n$ and each $|\pi |^j$-net is
finite, then the series $\sum_{j=1}^{\infty }y_jg_j$ converges. In
view of the inequalities above it converges to $f|_{W_q\setminus
W_{q-1}}$ by the norm of $L^{\infty }({\bf K}^m,\nu ,{\bf K})$ for
each $q$.  Denote $A_j$ for $f|_{W_q\setminus W_{q-1}}$ by $\mbox{
}_qA_j$. Then $\bigcup_{q=1}^{\infty }\mbox{ }_qA_j=:A_j$ is $\nu
$-measurable for each $j$. Since for each $x\in {\bf K}^m$ there
exists $q$ such that $x\in W_q\setminus W_{q-1}$, then
$\sum_{j=1}^{\infty }y_jg_j(x)=f(x)$ converges pointwise.
\par {\bf 7. Lemma.} {\it Suppose that $H$ is a compact non void
subset in ${\bf K}^n\times ({\bf K}\setminus \{ 0 \} )$ and $X_t :=
\{ y\in {\bf K}^n: |y-z|\le |x|^rt \mbox{ for every } (z,x)\in H \}
$ for $0\le t<\infty $, where $0<r\le 1$ is a constant, then
$c:=\inf \{ t: X_t\ne \emptyset \} <\infty $ and
$X_c=(\bigcap_{t>c}X_t)\ne \emptyset $. If $q\in X_c$, then $A_q\ne
\emptyset $, where $A_q := \{ z: \mbox{ there exists } (z,x)\in H
\mbox{ with } |q-z|=|x|^rc \} $.}
\par {\bf Proof.} Each set $X_t$ is compact, since $H$ is compact.
Consider the projection $\pi _{n+1}: {\bf K}^{n+1}\to \bf K$ such
that $\pi _{n+1}(\mbox{ }_1x,...,\mbox{ }_{n+1}x)= \mbox{ }_{n+1}x$,
where $\mbox{ }_jx\in \bf K$ for each $j=1,...,n+1$. Thus $\pi
_{n+1}(H)$ is compact as the continuous image of the compact set.
But $\pi _{n+1}(H)$ is contained in ${\bf K}\setminus \{ 0 \} $,
consequently, $\inf_{(z,x)\in H} |x|>0$. At the same time
$\sup_{(z,x)\in H} |z|<\infty $, hence $0\le \sup \{ |z|/|x|^r:
(z,x)\in H \} <\infty $. Therefore, $0\in X_t$ for each $t\ge \sup
\{ |z|/|x|^r: (z,x)\in H \} $. Then $X_c=\bigcap_{c<t<\infty }
X_t\ne \emptyset $, since $X_t\subset X_q$ for each $0<t\le q$ and
$X_t\ne \emptyset $ for each $t>c$. Put $R := \sup \{ |x|^r:
(z,x)\in H\mbox{ for some } z \} $. Consider $y, z\in X_c$, $\alpha
\in \bf K$ with $|\alpha |\le 1$, then $(a,q)\in H$ implies $|\alpha
y+(1-\alpha )z-a|\le \max (|\alpha | |y-a|, |1-\alpha | |z-a| )$,
but $\max (|\alpha |, |1-\alpha |)\le 1$, hence $|\alpha y+(1-\alpha
)z-a|\le |q|^rc$ and $\alpha y+(1-\alpha )z\in X_c$. Subjecting
${\bf K}^n$ to a translation we can suppose without loss of
generality that $0\in X_c$. If $z\in {\bf K}^n$ and $|z|=1$, then
$|\pi ^sz|>0$ for each $s\in \bf Z$. We have $\pi ^sz\in X_c$ if and
only if $|\pi ^sz-a|\le |x|^rc$ for each $(a,x)\in H$, but $|a|\le
|x|^rc$, since $0\in X_c$. Thus $\pi ^sz\in X_c$ if and only if
$|\pi ^sz|\le |x|^rc$ for each $(z,x)\in H$, which is equivalent to
$|\pi ^s|\le |x|^rc$ for each $(a,x)\in H$ and in its turn this is
equivalent to $|\pi ^s|\le c(\inf_{(z,x)\in H} |x|)^r$. Then there
exists $s_0\in \bf Z$ such that $|\pi ^s|\le  c(\inf_{(z,x)\in H}
|x|)^r$ for each $s\ge s_0$, $s\in \bf Z$. We have $H\cap \{ (z,x):
|z|=|x|^rc \} \ne \emptyset $, hence $A_q\ne \emptyset $, where
$q=0$ after translation.
\par {\bf 8. Theorem.} {\it If $S$ is a subset in ${\bf K}^m$ and
$f: S\to {\bf K}^n$ is a lipschitzian function with constants
$0<Lip_1(f) := C<\infty $ and $0< Lip_2(f) := r\le 1$: \par $(1)$
$|f(x)-f(y)\le C|x-y|^r$ for each $x, y\in S$, \\
then $f$ has a lipschitzian extension $g: {\bf K}^m\to \bf K$ such
that $Lip_1 (f)=Lip_1 (g)$ and $Lip_2(f)=Lip_2(g)$, where $n, m\in
\bf N$.}
\par {\bf Proof.} With the help of transformation $f\mapsto cf$,
where $0\ne c\in \bf K$ we can suppose that $0<Lip _1(f) =: b\le 1$.
Consider a class $\Psi $ of all lipschitzian extensions $f_j$ of $f$
on some subset $T_j$ of ${\bf K}^m$ having the same constants $C,
r$. Then $\Psi $ is partially ordered $(f_1,T_1)\preceq (f_2,T_2)$
if $T_1\subset T_2$ and $f_1|_{T_1}=f_2|_{T_1}$. Each linearly
ordered subset $\Phi $ in $\Psi $ has a maximal element $(h,T)$:
$h|_{T_j}=f_j$, $(f_j,T_j)\preceq (h,T)$ for each $(f_j,T_j)\in \Phi
$, where $T\supset \bigcup_jT_j$, since for each $T_j, T_k$ there is
the inequality $j\prec k$ if and only if $T_j\subset T_k$ and
$f_j|_{T_j}=f_k|_{T_j}$. In view of the Kuratowski-Zorn lemma
\cite{eng} there exists a maximal element $(g,T)$ in $\Psi $, $g:
T\to {\bf K}^n$, where $T\subset {\bf K}^m$. It is sufficient to
show, that if there exists $z\in {\bf K}^m\setminus T$, then there
exists $y\in {\bf K}^n$ such that $|y-g(x)|\le b |z-x|^r$ for every
$x\in T$, consequently, $g\cup \{ (z,y) \} \in \Psi $ and $(g,T)$
would not be maximal in $\Psi $. Thus we must prove, that
$\bigcap_{x\in T} B({\bf K}^n,g(x),b |x-z|^r)\ne \emptyset $. These
balls are compact, hence it is sufficient to prove that
$\bigcap_{x\in F} B({\bf K}^n,g(x),b |x-z|^r)\ne \emptyset $ for
each finite subset $F$ in $T$. Take $X_c$ from Lemma 7 and $q\in
X_c$, then $|q-g(x_i)|=|x_i-x|^rc$ for $i=1,...,k$, $g(x_i)\in A_q$,
where $1\le k\in \bf Z$, $(q,x)\in H$. We will show that $q\in X_b$.
If $0<c\le b$, then $q\in X_c\subset X_b$. If $c>b$, then
$|q-g(x_i)|>|x_i-x|^rb$ for each $x_i$, that is impossible by the
supposition of this theorem.
\par If $x$ is a limit point in $T$, then take a sequence $\{ x_n: n
\} $ such that $\lim_{n\to \infty }x_n=x$ and $q=\lim_{n\to \infty }
g(x_n)$, since $g$ is continuous on $T$ and $\{ g(x_n): n \} $ is
the Cauchy net in $\bf K$, but the latter uniform space is complete,
because $\bf K$ is a locally compact field. Then
$|q-g(y)|=\lim_{n\to \infty } |g(x_n)-g(y)|\le b \lim_{n\to \infty }
|x_n-y|^r=b|x-y|^r$ for each $y\in T$. Thus $T$ is the closed subset
in ${\bf K}^m$, hence ${\bf K}^m\setminus T$ is open. Suppose that
$v\in {\bf K}^m\setminus T$, then there exists $\delta :=\inf_{x\in
T} |v-x|>0$. Take $\delta \le R<\infty $, then $T\cap B({\bf
K},x_0,R)$ is compact, where $x_0\in T$ is a marked point.
Therefore, there exists $v_0\in T$ such that $|v-v_0|=\delta $.
\par We have that $g(T)$ is locally compact and closed in ${\bf K}^m$
and $|g(x)-g(y)|\le b|x-y|^r$ for each $x, y\in T$, consequently,
$B({\bf K}^m,g(x),b|x-y|^r)=B({\bf K}^m,g(y),b|x-y|^r)$ for each $x,
y\in T$. If $|y-v|>|x-v|$, then $|x-y|=|y-v|$, consequently, $B({\bf
K}^m,g(y),b|y-v|^r)=B({\bf K}^m,g(y),b|x-y|^r)=B({\bf
K}^m,g(x),b|x-y|^r)\supset B({\bf K}^m,g(x),b|x-v|^r)$, hence
\par $\bigcap_{x\in T} B({\bf K}^m,g(x),b|x-v|^r)\supset \bigcap_{x\in
T, |x-v|=|v_0-v|} B({\bf K}^m,g(x),b|x-v|^r)$. \\ On the other hand,
$|g(v_0)-g(x)|\le b|v_0-x|^r\le b\max
(|x-v|^r,|v_0-v|^r)=b|v_0-v|^r$ for $|x-v|=|v_0-v|$. Therefore,
$B({\bf K}^m,g(x),b|x-v|^r)=B({\bf K}^m,g(v_0),b|v_0-v|^r)$ and
inevitably $\bigcap_{x\in T} B({\bf K}^m,g(x),b|x-v|^r)\supset
B({\bf K}^m,g(v_0),b|v-v_0|^r)\ne \emptyset $, since the valuation
of $\bf K$ is non trivial and $\lim_{k\to\infty }|\pi |^k=0$.
\par {\bf 9. Theorem.} {\it Let $U$ be an open subset in $\bf K$
and let also $g: U\to \bf K$ be a locally lipschitzian function such
that for each $x_0\in U$ there exist constants $0<C<\infty $ and
$\delta >0$ and $0<r\le 1$ with \par $(1)$ $|g(x)-g(y)|\le C
|x-y|^r$ for each $\max (|x-x_0|, |y-x_0|)<\delta $, $x, y\in U$. \\
Then ${\bar {\Phi }}^1g(x;v;t)$ is continuous for $\mu ^3$-almost
all points in $U^{(1)}$ and $dg(x)/dx={\bar {\Phi }}^1g(x;1;0)$
exists and is continuous for $\mu $-almost all points of $U$.}
\par {\bf Proof.} The function $g$ is locally lipschitzian, hence
it is continuous. On the other hand, the Haar measure $\mu ^m$ on
${\bf K}^m$ is Radon, that is by the definition it satisfies the
following three conditions:
\par $(2)$ if $J$ is a compact subset of $\bf K$, then $\mu ^m
(J)<\infty $;
\par $(3)$ if $V$ is open in $\bf K$, then $\mu ^m(V)=
\sup \{ \mu ^m(J): J\mbox{ is compact }, J\subset V \} $; \par $(4)$
if $A$ is a $\mu ^m$-measurable subset, $A\subset \bf K$, then $\mu
^m (A) := \inf \{ \mu ^m(V): V \mbox{ is open }, A\subset V \} $.
\par In view of approximation Theorem 2.2.5 \cite{federer}
for each $\mu ^m $-measurable subset $A$ in ${\bf K}^m$ with $\mu
^m(A)<\infty $ and $\epsilon >0$ there exists a compact subset
$J\subset A$ such that $\mu ^m(A\setminus J)<\epsilon $. If $\cal E$
is a subset in $U^{(1)}$ of discontinuity of ${\bar {\Phi }}^1g$ and
$\cal D$ is a subset of discontinuity of $dg(x)/dx$ in $U$, then it
is sufficient to demonstrate, that $\mu ^3({\cal E}\cap
U_{R,\epsilon }^{(1)})=0$ and $\mu ({\cal D}\cap U_{R,\epsilon })=0$
for each $0<R<\infty $ and $\epsilon >0$, where $U_{R,\epsilon }$ is
a subset in $U\cap B({\bf K},0,R)$ such that $\mu (U\setminus
U_{R,\epsilon })<\epsilon $, since $\mu ^3((U^{(1)}\setminus
U_{R,\epsilon }^{(1)})\cap B({\bf K}^3,0,R))<3R^2\epsilon
+3R\epsilon ^2+\epsilon ^3$. For each compact set $U_{R,\epsilon }$
the covering $B({\bf K},x_0,\delta )$ with $\delta =\delta (x_0)>0$
has a finite subcovering, hence there exist $C=\sup_{x_0\in
U_{R,\epsilon }} C(x_0)<\infty $ and $0<r=\inf_{x_0\in U_{R,\epsilon
}} r(x_0)\le 1$ for which Inequality $(1)$ is satisfied for each $x,
y\in U_{R,\epsilon }$. Consider a restriction $g|_{U_{R,\epsilon
}}$, then by Theorem 8 it has a lipschitzian extension
$g_{R,\epsilon }$ on $\bf K$ with the same constants $0<C<\infty $
and $0<r\le 1$. Therefore, it is sufficient to prove this theorem
for a clopen compact subset $U$ in $\bf K$ which is supposed in the
proof below.
\par Since the mapping $x\mapsto x+vt$ is continuous by $(x,v,t)\in {\bf
K}^3$ and $\mu ^3$ is the Haar measure on ${\bf K}^3$ such that $\mu
^3$ has not any atoms, then ${\bar {\Phi
}}^1g(x;v;t)=[g(x+vt)-g(x)]/t$ is continuous for $\mu ^3$-almost all
points in $U^{(1)}$ if and only if $[g(x)-g(y)]/[x-y]$ is continuous
for $\mu ^2$-almost all points of ${\bf K}^2$.
\par Consider the following relation \par $V := \{ (x,S): S
\mbox{ is a compact clopen subset in } {\bf K}^m, x\in S \} $, where
$m\in \bf N$. Verify that $V$ is the $\mu ^m$ Vitaly relation.
Indeed, \par $(1)$ $V$ is a covering relation, that is a subset of $
\{ (x,S): x\in S\subset {\bf K}^m \} $;
\par $(2)$ $V\subset {\cal B}({\bf K}^m)$;
\par $(3)$ for each $0<R<\infty $ and each $y\in {\bf K}^m$
the ball $B({\bf K}^m,y,R)$ belongs to $V$, hence $\inf \{ diam (S):
(x,S)\in V \} =0$ for each $x\in {\bf K}^m$, consequently, $V$ is
fine at each point of ${\bf K}^m$;
\par $(4)$ ${\bf K}^m$ is locally compact separable and with a
countable base of its topology consisting of clopen balls. Thus if
$W\subset V$ and $Z\subset {\bf K}^m$ and $W$ is fine at each point
$z\in Z$, then $W(Z)$ has a countable disjoint subfamily covering
almost all of $Z$. Indeed, $W(Z)$ gives a base of topology inherited
from ${\bf K}^m$. This base is countable, hence $Z\subset
\bigcup_{j=1}^{\infty }U_j$, where each $U_j$ is a clopen compact in
${\bf K}^m$. Recall that a measure $\nu \in {\cal M}(X)$ is called
regular, if for each $A\subset X$ there exists a $\nu $-measurable
subset $G$ in $X$ such that $A\subset G$ and $\nu (A)=\nu (G)$. With
arbitrary measure $\nu $ one associates a regular measure by the
formula $\lambda (A) := \inf \{ \nu (G): A\subset G \mbox{ and }
G\mbox{ is } \nu \mbox{- measurable} \} $ (see Section 2.1.5 and the
Lusin's Theorem 2.3.5 \cite{federer}). Put $V_1=U_1$ and
$V_j=U_j\setminus \bigcup_{i<j}U_i$, then $V_i\cap V_j=\emptyset $
for each $i\ne j$. Since $\mu ^m$ is regular, then for each $V_j$
there exists a finite subfamily $W_{k,j}\in W$ such that $\mu
^m(V_j\setminus \bigcup_{k=1}^nW_{k,j})<\epsilon _j$, where
$n=n(\epsilon _j,j)\in \bf N$, $W_{k,j}\subset V_j$, $W_{k,j}\cap
W_{l,j}=\emptyset $ for each $k\ne l$. Choose $\epsilon _j=\epsilon
|\pi |^j$ and let $\epsilon >0$ tend to zero. Thus
$\bigcup_{k,j}W_{k,j}$ covers almost all of $Z$, since $0\le \mu
^1(Z\setminus \bigcup_{j,k}W_{k,j})\le \lim_{0<\epsilon \to 0}
\sum_j\mu ^m(V_j\setminus \bigcup_{k=1}^{n(\epsilon
_j,j)}W_{k,j})=0$.
\par Since the field $\bf K$ is locally compact, then there exists
a generator $|\pi |$ of the valuation group $\Gamma _{\bf K}$ such
that $|x|=|\pi |^{- \nu (x)}$ for each $x\in \bf K$, where
$|x|=mod_{\bf K}(x)$ is the multiplicative norm in $\bf K$, while
$\nu (x) = \nu _{\bf K}(x) \in \bf Z$ is called the valuation
function or valuation \cite{roo,sch1,weil}. For the first statement
of the theorem it is sufficient to prove, that the $\mu ^2$ measure
of the set \par ${\cal A} := \{ (x,y)\in U^2: \lim_{(x_1,y_1)\to
(x,y)} [g(x_1)-g(y_1)]/[x_1-y_1]$ \\ $\mbox{ either does not exist
or is not equal to } [g(x)-g(y)]/[x-y] \} $ \\  is zero, since $\mu
^2 \{ (x,x): x\in {\bf K} \} =0$. Since $g$ is the lipschitzian
function and
\par $[g(x_1)-g(y_1)]/[x_1-y_1] -
[g(x)-g(y)]/[x-y]=[(g(x_1)-g(y_1))((x-y)-(x_1-y_1)) +
((g(x_1)-g(y_1))-(g(x)-g(y)))(x_1-y_1)]/[(x_1-y_1)(x-y)]$, then
\par $(5)$ $|[g(x_1)-g(y_1)]/[x_1-y_1] - [g(x)-g(y)]/[x-y]|\le
C \max (|x_1-y_1|^r\max (|x_1-x|, |y_1-y|), |x_1-y_1| \max
(|x_1-x|^r, |y_1-y|^r))/[|x_1-y_1| |x-y|]$. \\
Under suitable affine mapping $q(x):=a(x-x_0)$ the image of $U$ is
contained in $B({\bf K},0,|\pi |)$, where $0\ne a\in \bf K$, $x_0\in
\bf K$, so without restriction of generality suppose that $U\subset
B({\bf K},0,|\pi |)$, since ${\bar {\Phi }}^1g$ is almost everywhere
continuous on $U^{(1)}$ if and only if ${\bar {\Phi }}^1g\circ
q^{-1}$ is such on $(q(U))^{(1)}$. Consider the sets
\par $(6)$ $A_{l,n,k}:=\{ (x,y)\in U^2: |x-y|=|\pi |^l,
\mbox{ there exists } (x_1,y_1)\in U^2 \mbox{ such that }$\\ $\max
(|x_1-x|, |y_1-y|)\le |\pi |^k, |[g(x_1)-g(y_1)]/(x_1-y_1)-
[g(x)-g(y)]/(x-y)|\ge |\pi |^n \} $, \\
where $l, n, k\in \bf N$. We have $\mu ^2((x,y)\in {\bf K}^2:
|x-y|=|\pi |^l, y\in B({\bf K},y_0,|\pi |^s))=(|\pi |^l-|\pi|^{l+1})
|\pi |^s$ for each $l, s\in \bf N$. If $|x-y|=|x_1-y_1|>\max
(|x_1-x|,|y_1-y|)$, then from $(5)$ it follows, that
$|[g(x_1)-g(y_1)]/[x_1-y_1] - [g(x)-g(y)]/[x-y]|\le C \max
(|x-y|^{(r-2)}\max (|x_1-x|, |y_1-y|),|x-y|^{-1}\max (|x_1-x|^r,
|y_1-y|^r))$.  Choose $k\in \bf N$ sufficiently large, $k\ge m_0$,
such that $C\max [|\pi |^{k+{(r-2)l}},|\pi |^{rk-l}]< |\pi |^n$,
then $\mu ^2(A_{l,n,k})=0$, since $\mu ^2(A_{l,n,k}\cap \{ y\in
B({\bf K},y_0, |\pi |^l) \} )=0$ for each $y_0\in \bf K$, where
$m_0=m_0(l)\in \bf N$. Let $l_0\in \bf N$ be a large number, take
$m_0=m_0(l_0)$ such that $l_0$ tends to the infinity if and only if
$m_0$ tends to the infinity, then $\mu ^2(\bigcup_{n=1}^{\infty
}\bigcup_{k, k\ge m_0} \bigcup_{l=1}^{\infty } A_{l,n,k})\le
\sum_{l=l_0}^{\infty } (|\pi |^{l+1}-|\pi |^{l+2})=|\pi |^{l_0+1}$,
since $B({\bf K},0,|\pi |)\setminus \{ 0 \} =\bigcup_{l=1}^{\infty }
\{ x\in {\bf K}: |x|=|\pi |^l \} $ and $\mu ( \{ 0 \} )=0$,
consequently, $\mu ^2(\bigcup_{n=1}^{\infty }\bigcap_{m=1}^{\infty
}\bigcup_{k, k\ge m} \bigcup_{l=1}^{\infty } A_{l,n,k})\le |\pi
|^{l_0+1}$, where $l_0$ is arbitrary large. Therefore, $\mu ^2(
\bigcup_{n=1}^{\infty } \bigcap_{m=1}^{\infty } \bigcup_{k, k\ge m}
\bigcup_{l=1}^{\infty } A_{l,n,k})=0$, consequently, $\mu ^2({\cal
A})=0$, since due to $(6)$ \par ${\cal A}\subset \{ (x,y)\in U^2:
\mbox{ there exists a sequence } (x^m,y^m)\in U^2 \mbox{ such that }
$ \\  $\lim_{m\to \infty } (x^m,y^m)=(x,y), {\overline {\lim
}}_{m\to \infty } |[g(x^m)-g(y^m)]/[x^m-y^m] - [g(x)-g(y)]/[x-y]|\ge
|\pi |^n \mbox{ for some } n\in {\bf N} \} \subset (
\bigcup_{n=1}^{\infty } \bigcap_{m=1}^{\infty } \bigcup_{k, k\ge m}
\bigcup_{l=1}^{\infty } A_{l,n,k})$. \\ Thus $[g(x)-g(y)]/(x-y)$ is
$\mu ^2$-almost everywhere continuous on $U^2$.
\par Now prove the second statement, for this mention that the set
$\cal A$ is symmetric relative to the transposition $(x,y)\mapsto
(y,x)$. We have that $[g(x)-g(y)]/(x-y)$ is continuous for $\mu
^2$-almost all $(x,y)\in U^2$, hence on everywhere dense subset in
$U^2$. It is sufficient to show that the $\mu $ measure of the set
${\cal C} := \{ x\in U: \lim_{(x_1,y_1)\to (x,x)}
[g(x_1)-g(y_1)]/[x_1-y_1]$ \\
$\mbox{ either does not exist or is not equal to } \lim_{y\to x}
[g(x)-g(y)]/[x-y]$\\  $\mbox{ or the latter limit does not exist }
\} $ \\ is zero. For this consider the sets
\par $(7)$ $E_{l,n,k}:=\{ x\in U:
\mbox{ there exist } y\in U \mbox{ and }
(x_1,y_1)\in U^2 \mbox{ such that }$ \\
$|x-y|=|\pi |^l, \max (|x_1-x|, |y_1-y|)\le |\pi |^k,
|[g(x_1)-g(y_1)]/(x_1-y_1)-
[g(x)-g(y)]/(x-y)|\ge |\pi |^n \} $, \\
where $l, n, k\in \bf N$. There exists $m_0\in \bf N$ such that for
each $k\ge m_0$ there is satisfied the inequality $C\max [|\pi
|^{k+{(r-2)l}},|\pi |^{rk-l}]< |\pi |^n$, hence $\mu (E_{l,n,k})=0$
for such $k$, since $\mu ^1(E_{l,n,k}\cap \{ y\in B({\bf K},y_0,
|\pi |^l) \} )=0$ for each $y_0\in \bf K$. For $l_0\in \bf N$ take
$m_0=m_0(l_0)$ such that $l_0$ tends to the infinity if and only if
$m_0$ tends to the infinity, then
\par $\mu (\bigcup_{n=1}^{\infty }\bigcup_{k, k\ge m} \bigcup_{l,
l\ge l_0} E_{l,n,k})\le \sum_{l=l_0}^{\infty } (|\pi |^l-|\pi
|^{l+1})=|\pi |^{l_0}$, consequently, $\mu ^1(\bigcup_{n=1}^{\infty
}\bigcap_{m=1}^{\infty }\bigcup_{k, k\ge m} \bigcap_{s=1}^{\infty
}\bigcup_{l, l\ge s} E_{l,n,k})\le |\pi |^{l_0}$, where $l_0$ is
arbitrary large. Therefore, $\mu ( \bigcup_{n=1}^{\infty }
\bigcap_{m=1}^{\infty } \bigcup_{k, k\ge m} \bigcap_{s=1}^{\infty
}\bigcup_{l, l\ge s} E_{l,n,k})=0$, consequently, $\mu ({\cal
C})=0$, since in view of $(7)$ \par ${\cal C}\subset \{ x\in U:
\mbox{ there
exists a sequence } (z^m,x^m,y^m)\in U^3 \mbox{ such that } $ \\
$\lim_{m\to \infty } (z^m,x^m,y^m)=(x,x,y), {\overline {\lim
}}_{m\to \infty } |[g(x^m)-g(y^m)]/[x^m-y^m] -
[g(x)-g(z^m)]/[x-z^m]|\ge |\pi |^n \mbox{ for some } n\in {\bf N} \}
\subset ( \bigcup_{n=1}^{\infty } \bigcap_{m=1}^{\infty }
\bigcup_{k, k\ge m} \bigcap_{s=1}^{\infty }\bigcup_{l, l\ge s}
E_{l,n,k})$. Therefore, $dg(x)/dx$ exists and is $\mu $-almost
everywhere continuous on $U$.
\par {\bf 10. Lemma.} {\it If $S$ is a $\mu ^m\otimes \mu ^k$
measurable subset of ${\bf K}^m\times {\bf K}^k$, $\epsilon >0$ and
$\delta >0$ and $T:= \{ x\in {\bf K}^m: \mu ^k(\{ z: (x,z)\in S,
|z|\le R \} )\le \epsilon R^k \} $, whenever $0<R<\delta $, then $T$
is $\mu ^m$ measurable. }
\par {\bf Proof.} For each $0<R<\infty $ the set $S_R:= S\cap \{
(x,z): |z|\le R \} $ is $\mu ^m\otimes \mu ^k$ measurable and by the
Fubini theorem $\mu ^k ( \{ z: (x,z)\in S_R \} )$ is the $\mu ^m$
measurable function of the variable $x$, hence $T$ is $\mu ^m$
measurable, since $T= \{ x\in {\bf K}^m: \mu ^k(\{ z: (x,z)\in S_R
\} )\le \epsilon R^k \} $.
\par {\bf 11. Lemma.} {\it If a function $\phi : {\bf K}^m\times
{\bf K}^k \to \bf R$ is $\mu ^m\otimes \mu ^k$ measurable, then $ap
{\overline {\lim}}_{z\to 0} \phi (x,z)$ and $ap {\underline
{\lim}}_{z\to 0} \phi (x,z)$ are $\mu ^m$ measurable functions of
the variable $x$.}
\par {\bf Proof.} For each $c\in \bf R$ applying Lemma 10 to the
sets $ \{ (x,z): \phi (x,z)>c \} $ and $\{ (x,z): \phi (x,z)<c \} $
we get the statement of this lemma.
\par {\bf 12. Lemma.} {\it If $u: {\bf K}^n\to {\bf K}^m$ is a $\bf
K$-linear epimorphism and $A$ is a $\mu ^m$ measurable set, then
$u^{-1}(A)$ is $\mu ^n$ measurable.}
\par {\bf Proof.} This follows from the fact that there exists
a $\bf K$-linear isomorphism $v: {\bf K}^n\to {\bf K}^m\times {\bf
K}^{n-m}$ such that $v \circ u^{-1}(A)=A\times {\bf K}^{n-m}$.
\par {\bf 13. Lemma.} {\it If a function $g: {\bf K}^m\to \bf K$
is $\mu ^m$ measurable and $1\le k\le m$, then the $\mu ^k$
approximate limit is $ap \lim_{z\to 0} g(\mbox{ }_1x+\mbox{
}_1z,...,\mbox{ }_kx+\mbox{ }_kz,\mbox{ }_{k+1}x,...,\mbox{
}_mx)=g(x)$ for $\mu ^m$ almost all $x$.}
\par {\bf Proof.} There exists the isomorphism $v: {\bf K}^m
\to {\bf K}^k\times {\bf K}^{m-k}$. Put $h(x,y,z)=(x+z,y)$, where
$h: {\bf K}^k\times {\bf K}^{m-k}\times {\bf K}^k\to {\bf K}^k\times
{\bf K}^{m-k}$. In view of Lemma 12 the composite function $g\circ
h$ is $\mu ^k\times \mu ^{m-k}\times \mu ^k$ measurable. From Lemma
11 it follows, that $A:= \{ (x,y): ap \lim_{z\to 0}g(x+z,y)=g(x,y)
\} $ is $\mu ^k\otimes \mu ^{m-k}$ measurable. \par Theorem 2.9.13
\cite{federer} states that if a function $f$ maps $\nu $ almost all
of $X$ into $Y$, where $(X,\rho _X)$ is a metric space, $(Y,\rho
_Y)$ is a separable metric space, then $f$ is $\nu $ measurable if
and only if $f$ is approximately $(\nu ,V)$ continuous at almost all
points of $X$. \par In accordance with the latter theorem and the
Fubini theorem the function $g(x,y)$ is $\mu ^k$ measurable by $x$
for $\mu ^{m-k}$ almost all $y$ and $\mu ^k( \{ x: (x,y)\notin A \}
)=0$ and the complement of $A$ has $\mu ^k$ measure zero.
\par {\bf 14. Corollary.} {\it If $A$ is $\mu ^m$ measurable set in
${\bf K}^m$ and $1\le k\le m$, then for $\mu ^m$ almost all $y\in A$
the set ${\bf K}^k\cap \{ x: (\mbox{ }_1x,...,\mbox{ }_kx,\mbox{
}_{k+1}y,..., \mbox{ }_my)\notin A \} $ has zero $\mu ^k$ density at
$(\mbox{ }_1y,...,\mbox{ }_ky)$.}
\par {\bf Proof.} This follows from Lemma 13 for the
characteristic function $g=ch_A$, where $ch_A(y)=1$ for each $y\in
A$, $ch_A(y)=0$ for each $y\in {\bf K}^m\setminus A$.
\par {\bf 15. Theorem.} {\it If $f: {\bf K}^m\to {\bf K}^n$ is
$\mu ^m$ measurable, then $A_i:=dom \mbox{ }_{ap}D_if$ is a $\mu ^m$
measurable set, $V_i:= dom \mbox{ }_{ap} {\bar {\Phi }}^1f(x;e_i;t)$
is a $\mu ^{m+1}$ measurable set such that $\mu ^{m+1}({\bf
K}^{m+1}\setminus V_i)=0$, moreover, $\mbox{ }_{ap}D_if$ and $\mbox{
}_{ap} {\bar {\Phi }}^1f(x;e_i;t)$ are $\mu ^m|_{A_i}$ and $\mu
^{m+1}|_{V_i}$ measurable functions respectively,\par $(1)$ $\mbox{
}_{ap} {\bar {\Phi }}^1f(x;v;t)=v_1\mbox{ }_{ap} {\bar {\Phi
}}^1f(x+e_2v_2+...+e_mv_m;e_1;v_1t)+ v_2\mbox{ }_{ap} {\bar {\Phi
}}^1f(x+e_3v_2+...+e_mv_m;e_2;v_2t)+...+ v_m\mbox{ }_{ap} {\bar
{\Phi }}^1f(x;e_m;v_mt)$ \\
for $\mu ^{m+1}$ almost all points $(x,t)$ in
$V:=\bigcap_{i=1}^mV_i$ and each $v=v_1e_1+...+v_me_m\in {\bf K}^m$,
\par $(2)$ $\mbox{ }_{ap}Df(x).v=\sum_{i=1}^mv_i\mbox{
}_{ap}D_if(x))$ \\
for $\mu ^m$ almost all points $x$ in $A := \bigcap_{i=1}^mA_i$ and
each $v\in {\bf K}^m$.}
\par {\bf Proof.} Since $f=(f_1,...,f_n)$, where
$f_j: {\bf K}^m\to \bf K$, then $dom \mbox{ }_{ap} {\bar {\Phi
}}^1f(x;e_i;t)=\bigcap_{j=1}^n dom \mbox{ }_{ap} {\bar {\Phi
}}^1f_j(x;e_i;t)$ and $dom \mbox{ }_{ap}D_if=\bigcap_{j=1}^n dom
\mbox{ }_{ap}D_if_j$ for each $i=1,...,m$, consequently, it is
sufficient to prove this theorem for $n=1$. Thus suppose that $n=1$.
\par In accordance with Theorem 2.9.13 \cite{federer} (see its
formulation in section 13) the function ${\bar {\Phi }}^1f(x;v;t)$
is approximately continuous on ${\bf K}^m\times {\bf K}^m\times
({\bf K}\setminus \{ 0 \} )$, since $[f(x+vt)-f(x)]/t$ is $\mu
^{2m+1}$ measurable on ${\bf K}^m\times {\bf K}^m\times ({\bf
K}\setminus \{ 0 \} )$ and inevitably $\mu ^{m+1}({\bf
K}^{m+1}\setminus V_i)=0$, since $\mu (\{ 0 \} )=0$. Therefore,
Formula $(1)$ is satisfied for $\mu ^{m+1}$ almost all points
$(x;t)$ in $V$, where $\mu ^{m+1}({\bf K}^{m+1}\setminus V)=0$. So
it remains to spread this for $t=0$, but $\mbox{ }_{ap}{\bar {\Phi
}}^1f(x;v;0)=ap \lim_{t\to 0} {\bar {\Phi }}^1f(x;v;t)=\mbox{
}_{ap}Df(x).v$ and the proof of this theorem reduces to the proof of
its statement relative to $\mbox{ }_{ap}Df(x)$. \par The set $A_i$
is $\mu ^m$ measurable if and only if its complement ${\bf
K}^m\setminus A_i$ is $\mu ^m$ measurable. But ${\bf K}^m\setminus
A_i=\bigcup_{s=1}^{\infty }\bigcap_{l, l\ge s} \bigcup_{u=1}^{\infty
}\bigcap_{k, k\ge u}\bigcup_{q=1}^{\infty }E_{l,k,q}$, where
$E_{l,k,q}:= \{ x\in {\bf K}^m: \mbox{ there exist } t_1, t_2\in
{\bf K} \mbox{ such that } \max (|t_1|,|t_2|)=|\pi |^k,
|t_1-t_2|=|\pi |^l, |\mbox{ }_{ap}{\bar {\Phi }}^1f(x;e_i;t_1) -
\mbox{ }_{ap}{\bar {\Phi }}^1f(x;e_i;t_2)|\ge |\pi |^q \mbox{ or }
|\mbox{ }_{ap}{\bar {\Phi }}^1f(x;e_i;t_1)|\ge |\pi |^{-q} \} $.
Each set $E_{l,k,q}$ is $\mu ^m$ measurable, since $\mu ^{m+1}({\bf
K}^{m+1}\setminus V_i)=0$ and due to Lemma 10, hence $A_i$ is $\mu
^m$ measurable for each $i=1,...,m$.
\par Put $T_{R,i,j}(x):= \{ t\in {\bf K}: |t|<R, |f(x+te_i)-f(x)-t
\mbox{ }_{ap}D_if(x)|>|t| |\pi |^j \} $ and $B_{i,j,q}:= A_i\cap \{
x\in {\bf K}^m: \mu (T_{R,i,j}(x))\le R|\pi |^j \forall 0<R<|\pi |^q
\} $, whenever $x\in {\bf K}^m$, $0<R<\infty $, $i, j, q\in \bf N$.
For $i>1$ denote $Z_{R,i,j,q}(x) := \{ z\in {\bf K}^{i-1}: |z|<R,
\mbox{ either } x+\mbox{ }_1ze_1+...+\mbox{ }_{i-1}ze_{i-1}\notin
B_{i,j,q} \mbox{ or } |\mbox{ }_{ap}D_if(x+\mbox{ }_1ze_1+...+\mbox{
}_{i-1}ze_{i-1}) - \mbox{ }_{ap}D_if(x)|>R|\pi |^j \} $, also
introduce the sets $C_{i,j,q,k}:=B_{i,j,q}\cap \{ x: \mu ^{i-1}
(Z_{R,i,j,q}(x))\le R^{i-1}|\pi |^j \forall 0<R<|\pi |^k \} $ for
$k\in \bf N$ and $i>1$, in particular, $C_{1,j,q,k}:=B_{1,j,q}$ for
each $k$. In view of Lemmas 10 and 13 $B_{i,j,q}$ and $C_{i,j,q,k}$
are  $\mu ^m$ measurable and $A_i=\bigcup_{q=1}^{\infty }B_{i,j,q}$
for each $(i,j)$, moreover, $\mu ^m(B_{i,j,q}\setminus
\bigcup_{k=1}^{\infty }C_{i,j,q,k})=0$ for each $(i,j,q)$, also
$B_{i,j,q}\subset B_{i,j,q+1}$ and $C_{i,j,q,k}\subset
C_{i,j,q,k+1}$. For a subset $S$ in $\bigcap_{i=1}^mA_i$ with $\mu
^m(S)<\infty $ and every $\epsilon >0$ there exist sequences $\{
q_j: j\in {\bf N} \} $ and $ \{ k_j: j\in {\bf N} \} $ of natural
numbers such that $\mu ^m (S\setminus B_{i,j,q_j})<\epsilon |\pi
|^j$ and $\mu ^m(S\cap B_{i,j,q_j}\setminus
C_{i,j,q_j,k_j})<\epsilon |\pi |^j $ for $i=1,...,m$ and $j\in \bf
N$, consequently, $\mu ^m(S\setminus G)< 2|\pi |(1-|\pi |)^{-1}
m\epsilon $, where $G:=\bigcap_{i=1}^m\bigcap_{j=1}^{\infty
}C_{i,j,q_j,k_j}$.
\par We will demonstrate that $f$ is uniformly approximately
differentiable at the points of $G$. Consider $x\in G$, $j\in \bf N$
and $0<R<\min (|\pi |^{q_j},|\pi |^{k_j})$, $S_i:= \{ v\in {\bf
K}^m: |v|\le R \mbox{ and either } i>1 \mbox{ and }
(v_1,...,v_{i-1})\in Z_{R,i,j,q_j}(x) \mbox{ or } v_i\in
T_{R,i,j}(x+v_1e_1+...+v_{i-1}e_{i-1}) \} $, where $i=1,...,m$.
Since $\mu ^{i-1} (Z_{R,i,j,q_j}(x))\le R^{i-1}|\pi |^j$ and since
$\mu (T_{R,i,j}(x))\le R|\pi |^j$ for
$z=x+v_1e_1+...+v_{i-1}e_{i-1}\in B_{i,j,q_j}$, then $\mu ^m(S_i)\le
R^{i-1}|\pi |^jR^{m-i+1}+R|\pi |^jR^{m-1}=2R^m|\pi |^j.$ If $v\in
B({\bf K}^m,0,R)\setminus S_i$, then $|f(x+v_1e_1+...+v_ie_i) -
f(x+v_1e_1+...+v_{i-1}e_{i-1})-v_i\mbox{ }_{ap}D_if(x)|\le \max
(|v_i||\pi |^j, |v_i\mbox{ }_{ap}D_if(x+v_1e_1+...+v_{i-1}e_{i-1}) -
v_i\mbox{ }_{ap}D_if(x)|)\le |v_i| |\pi |^j$. Putting
$S:=\bigcup_{i=1}^mS_i$ we get $\mu ^m(S)R^{-m}\le 2m|\pi |^j$ and
the inclusion $v\in B({\bf K}^m,0,R)\setminus S$ implies that
$|f(x+v)-f(x)-\sum_{i=1}^mv_i\mbox{ }_{ap}D_if(x)|\le |\pi |^j
\max_{i=1}^m |v_i|=|\pi |^j|v|$.
\par {\bf 16. Lemma.} {\it Let $S\subset A\subset {\bf K}^m$,
$f: A\to {\bf K}^n$, $0<R<\infty $, $0<C<\infty $ and $z\in S$
implies $B({\bf K}^m,z,R)\subset A$ with $|f(x)-f(z)|\le C|x-z|$ for
each $x\in B({\bf K}^m,z,R)$. Suppose also that $y\in S$, ${\bf
K}^m\setminus S$ has $\mu ^m$ density $0$ at each $z\in U_y$ for
some (open) neighborhood $U_y$ of $y$ in ${\bf K}^m$ and $f$ is
approximately differentiable at $y$. Then $f$ is differentiable at
$y$ such that ${\bar {\Phi }}^1f(y;v;t)$ is continuous on $(U_y\cap
A)^{(1)}$ for some neighborhood $U_y$ of $y$ in ${\bf K}^m$.}
\par {\bf Proof.} Suppose that $L=\mbox{ }_{ap}Df(y)$, $0<\epsilon
<1$, $0<\delta \le R$ and put $W:=S\cap \{ z: |f(z)-f(y)-L(z-y)|\le
\epsilon |z-y| \} $ and $\mu ^m(B({\bf K}^m,y,q)\setminus
W)<\epsilon ^mR^m$ for each $0<q<\delta $, $q\in \Gamma _{\bf K}$.
For $x\in B({\bf K}^m,y,\delta )$ take $q=|x-y|$ and mention that
$B({\bf K}^m,x,q)=B({\bf K}^m,y,q)$, also $W\cap B({\bf K}^m,x,q)\ne
\emptyset $. Choose $z\in B({\bf K}^m,x,\epsilon q)\cap W$, hence
$x\in B({\bf K}^m,z,\epsilon q)\subset B({\bf K}^m,z,R)$,
consequently,
\par $(1)$ $|f(x)-f(y)-L(x-y)|\le \max (|f(z)-f(y)-L(z-y)|,
|f(x)-f(z)|,|L(z-x)|)\le \max (\epsilon |z-y|,C|x-z|, \| L \|
|x-z|)\le \epsilon \max (|x-y|, Cq, \| L \| q)=\epsilon |x-y|\max
(1,C, \| L \| )$, since $|z-y|\le \max (|z-x|,|x-y|)\le |x-y|=q$.
From the arbitrariness of $0<\epsilon <1$ it follows, that there
exists $L=Df(y)$. From Inequality $(1)$ we have
\par $(2)$ $|{\bar {\Phi }}^1f(y;v;t) - Lv |\le \epsilon \max
(1, C, \| L \| )$ \\
for all $v=x-y$, $x\in B({\bf K}^m,y,\delta )$, consequently,
\par $(3)$ $\lim_{t\to 0} {\bar {\Phi }}^1f(y;v;t)=Lv$ uniformly
by $x\in B({\bf K}^m,y,\delta )$,\\ since $\epsilon >0$ is arbitrary
small. From the existence of $Df(y)$ it follows that $f$ is
continuous at $y$. The condition ${\bf K}^m\setminus S$ has $\mu ^m$
density $0$ at each $z\in U_y$ implies that $S$ is everywhere dense
in $U_y$ for some neighborhood $U_y$ of $y$. Then
$|f(y+xt_1)-f(y+vt_2)|\le \max ( |f(y+xt_1)-f(z)|,
|f(z)-f(y+vt_2)|)\le C\max (|z-y-xt_1|, |z-y-vt_2|)$ for each $z\in
S$, where $y+xt_1, y+vt_2\in A$, $B({\bf K}^m,z,R)\subset A$ and
$\max (|z-y-xt_1|, |z-y-vt_2|)$ can be chosen equal to
$|xt_1-vt_2|$, since $\Gamma _{\bf K}$ is discrete in $(0,\infty )$
and $S$ is dense in $U_y$. Thus ${\bar {\Phi }}^1f(y;v;t)$ is
continuous on $U_y^{(1)}\cap (A^{(1)}\setminus \{ (y;v;t): {\bf
K}^m\ni v\ne 0, {\bf K}\ni t\ne 0 \} )$. Together with $(3)$ this
gives the continuity of ${\bar {\Phi }}^1f(y;v;t)$ on $(U_y\cap
A)^{(1)}$.
\par {\bf 17. Theorem.} {\it If $f: {\bf K}^m\to {\bf K}^n$ is
a locally lipschitzian function such that for each $x_0\in {\bf
K}^m$ there exist constants $0<C<\infty $ and $\delta >0$ and
$0<r\le 1$ with
\par $(1)$ $|g(x)-g(y)|\le C |x-y|^r$ for each $\max (|x-x_0|,
|y-x_0|)<\delta $, $x, y\in {\bf K}^m$.\\
Then $f$ is differentiable at $\mu ^m$ almost all points of ${\bf
K}^m$ on a subset $G$ and ${\bar {\Phi }}^1f(x;v;t)$ is continuous
at $\mu ^{2m+1}$ almost all points of ${\bf K}^{2m+1}$ on a subset
$V_0$ such that for each $\epsilon >0$ there exists a closed subset
$G_{\epsilon }$ in $G$ with $\mu ^m(G\setminus G_{\epsilon
})<\epsilon $ and $Df$ is continuous on $G_{\epsilon }$ and
$G_{\epsilon }^{(1)}\subset V_0$.}
\par {\bf Proof.} Let the sets $A_i$ and $V_i$ be the same as in
Theorem 15. For $x\in {\bf K}^m$ consider the map $f_x(h):= f(\mbox{
}_1x,...,\mbox{ }_{i-1}x,h,\mbox{ }_{i+1}x,...,\mbox{ }_mx)$ for any
$h\in \bf K$, so it is locally lipschitzian. In accordance with
Theorem 9 ${\bar {\Phi }}^1f_x(h;w;t)$ is continuous for $\mu ^3$
almost every $(h;w;t)\in {\bf K}^3$ and $df_x(h)/dh$ is $\mu $
almost everywhere continuous by $h$ on ${\bf K}$ for each marked
point $(\mbox{ }_1x,...,\mbox{ }_{i-1}x,\mbox{ }_{i+1}x,...,\mbox{
}_mx)\in {\bf K}^{m-1}$. Since $A_i$ is $\mu ^m$ measurable and
$V_i$ is $\mu ^{2m+1}$ measurable, then $\mu ^m({\bf K}^m\setminus
A_i)=0$ and $\mu ^{2m+1}({\bf K}^{2m+1}\setminus V_i)=0$ and
inevitably $\mu ^m({\bf K}^m\setminus G)=0$ and $\mu ^{2m+1}({\bf
K}^{2m+1}\setminus V)=0$, where $G := \bigcap_{i=1}^mA_i$ and
$V=\bigcap_{i=1}^mV_i$. Thus $f|_G$ is differentiable, ${\bar {\Phi
}}^1f(x;v;t)|_V$ is $\mu ^{2m+1}$ measurable. On the other hand,
$({\bf K}^{2m+1}\setminus {\bf K}^{2m}\times \{ 0 \} ) \subset V$,
since $f$ is continuous on ${\bf K}^m$ and so ${\bar {\Phi
}}^1f(x;v;t)$ is continuous on $({\bf K}^{2m+1}\setminus {\bf
K}^{2m}\times \{ 0 \} )$, consequently, ${\bar {\Phi }}^1f(x;v;t)$
is continuous by each triple $(\mbox{ }_jx;\mbox{ }_jv;t)$ and $\mu
^{2m+1}$ measurable on $G^{(1)}$, $G^{(1)}\subset V$, $\mu ^{2m+1}
({\bf K}^{2m+1}\setminus G^{(1)})=0$, since $\mu ^m({\bf
K}^m\setminus G)=0$.
\par The Lusin Theorem 2.3.5 \cite{federer} asserts:
if $\phi $ is a Borel regular nonnegative measure over a metric
space $X$ (or a Radon measure over a locally compact Hausdorff space
$X$), if $f$ is a $\phi $ measurable function with values in a
separable metric space $Y$, $A$ is a $\phi $ measurable set for
which $\phi (A)<\infty $, and $\epsilon >0$, then $A$ contains a
closed (compact) set $C$ such that $\phi (A\setminus C)<\epsilon $
and $f|_C$ is continuous.
\par In view of the Lusin theorem for each $\epsilon
>0$ and each $0<R<\infty $ and $\xi \in G$ there exists a
compact subset $E\subset (G\cap B({\bf K}^m,\xi ,R))$ such that $\mu
^m((G\cap B({\bf K}^m,\xi ,R)) \setminus E)<\epsilon $ and $Df|_E$
is continuous, hence ${\bar {\Phi }}^1f(x;v;t)|_{E^{(1)}}$ is
continuous. Therefore, the restriction of $f$ on $E$ is lipschitzian
with $Lip_2(f|_E)=1$. In view of Theorem 8 there exists an extension
$g_E$ on ${\bf K}^m$ of $f|_E$ with $Lip_1(g_E)=Lip_1(f|_E)$ and
$Lip_2(g_E)=1$. \par  Take in Lemma 16 $S=E$ and $A={\bf K}^m$.
Therefore, $g_E$ is differentiable at $\mu ^m$ almost all points of
${\bf K}^m$ and ${\bar {\Phi }}^1g_E(x;v;t)$ is continuous at $\mu
^{2m+1}$ almost all points of ${\bf K}^{2m+1}$. Since $\xi \in G$,
$0<R<\infty $ and $\epsilon >0$ are arbitrary, we can take a
disjoint covering $B({\bf K}^m,\xi _j,R_j)$, $R_j\ge 1$, of $G$ and
$E_{i,j}\subset [(G\cap B({\bf K}^m,\xi _j,R_j))\setminus
\bigcup_{l, l<i}E_{l,j}]$ such that $\mu ^m((G\cap B({\bf K}^m,\xi
_j,R_j))\setminus (\bigcup_{l=1}^i E_{l,j}))<|\pi |^{i+j+k}$, where
$k\in \bf N$ is some large fixed number, $i, j\in \bf N$,
$E_{0,j}:=\emptyset $. Then consider $G\setminus
\bigcup_{j=1}^s\bigcup_{i=1}^sE_{i,j}$ and continue this
construction by induction. The family of restrictions $f|_{E_{i,j}}$
generates the function $f|_{G_1}$, where
$G_1=\bigcup_{i,j=1}^{\infty }E_{i,j}$ and $\mu ^m(G\setminus
G_1)=0$. Since $E_{i,j}\subset [B({\bf K}^m,\xi _j,R_j))\setminus
\bigcup_{l, l<i}E_{l,j}]$, each $E_{l,j}$ is compact and each
$[B({\bf K}^m,\xi _j,R_j))\setminus \bigcup_{l, l<i}E_{l,j}]$ is
open in ${\bf K}^m$ and $B({\bf K}^m,\xi _j,R_j)\cap B({\bf K}^m,\xi
_q,R_q)=\emptyset $ for each $q\ne j$, then take $G_{\epsilon
}:=\bigcup_{j=1}^{\infty }\bigcup_{i=1}^{n(j)}E_{i,j}$, where
$n(j)\in \bf N$ is a sequence such that $\sum_{j=1}^{\infty }|\pi
|^{n(j)+j+k}<\epsilon $. Each $n(j)$ is finite and
$\bigcup_{i=1}^{n(j)}E_{i,j}$ is closed, hence
$\bigcup_{j=1}^{\infty }\bigcup_{i=1}^{n(j)}E_{i,j}$ is closed in
${\bf K}^m$ and inevitably $Df$ is continuous on $G_{\epsilon }$
such that $G_{\epsilon }^{(1)}\subset V_0$.
\par {\bf 18. Lemma.} {\it  If $A$ is a $\mu ^m$-measurable subset
in ${\bf K}^m$ and $f: A\to {\bf K}^n$ is locally lipschitzian, then
$f$ has ${\mu }^{2m+1}$-everywhere in $A^{(1)}$ an approximate
partial difference quotient ${\bar {\Phi }}^1f(x;v;t)$ and ${\bar
{\Phi }}^1f(x;e_i;t)$ for $\mu ^{m+1}$ almost all points $(x;e_i;t)$
of $(A\times \{ e_i \} \times {\bf K})\cap A^{(1)}$ and $\mu
^m$-almost everywhere on $A$ an approximate differential $Df(x)$ and
approximate partial differentials $D_jf(x)=\partial f(x)/\partial
\mbox{ }_jx$.}
\par {\bf Proof.} In accordance with Theorem 8 and the proof of Theorem
9 for each $\epsilon >0$ and $0<R<\infty $ the function
$f|_{A_{R,\epsilon }}$ has a lipschitzian extension $g: {\bf K}^m\to
{\bf K}^n$ such that $g$ has the same lipschitzian constants
$0<C<\infty $ and $0<r\le 1$, where $A_{R,\epsilon }$ is a compact
subset in $A\cap B({\bf K}^m,0,R)$ such that $\mu ^m(A\cap B({\bf
K}^m,0,R)\setminus A_{R,\epsilon })<\epsilon $. From Theorem 17 it
follows that $g$ has ${\bar {\Phi }}^1g$ and $Dg$ at $\mu ^{2m+1}$
and $\mu ^m$ almost all points of $A^{(1)}$ and $A$ respectively.
Recall that one say that $B$ is a $\phi $ hull of $A$ if and only if
$A\subset B\subset X$, $B$ is $\phi $ measurable and $\phi (T\cap
A)=\phi (T\cap B)$ for every $\phi $ measurable subset $T$, where
$\phi $ is a measure over $X$. Theorem 2.9.11 \cite{federer} states
that if $A\subset X$ and $P=\{ x: (V)\lim_{S\to x} \phi (S\cap
A)/\phi (S)=1 \} $, $Q= \{ x: (V)\lim_{S\to x} \phi (S\setminus
A)/\phi (S)=0 \} $, then $P$ and $Q$ are $\phi $ measurable, $\phi
(A\setminus P)=0$, $A\cup P$ is a $\phi $ hull of $A$, $\phi
(Q\setminus A)=0$, $X\setminus (A\cap Q)$ is a $\phi $ hull of
$X\setminus A$. Moreover, $\phi $ measurability of $A$ is equivalent
to each of the two conditions: $\phi (P\setminus A)=0$, $\phi
(A\setminus Q)=0$.
\par In view of this theorem ${\bf K}^{2m+1}\setminus
A^{(1)}$ and ${\bf K}^m\setminus A$ have zero densities at $\mu
^{2m+1}$ and $\mu ^m$ almost all points of $A^{(1)}$ and $A$
respectively. At points where both conditions hold there are $\mbox{
}_{ap}{\bar {\Phi }}^1g$ and $\mbox{ }_{ap}Dg$. By Corollary 14 we
have $D_ig(x)=\mbox{ }_{ap}D_ig(x)$ and ${\bar {\Phi }}^1g(x;e_i;t)=
\mbox{ }_{ap}{\bar {\Phi }}^1g(x;e_i;t)$ for $\mu ^m$ almost all
$x\in A$ and $\mu ^{m+1}$ almost all points $(x;e_i;t)$ of $(A\times
\{ e_i \} \times {\bf K})\cap A^{(1)}$ correspondingly. From
combinatorial Formulas $15(1,2)$ we get, that $f$ has ${\mu
}^{2m+1}$-everywhere in $A^{(1)}$ an approximate partial difference
quotient $\mbox{ }_{ap}{\bar {\Phi }}^1f(x;v;t)$ and $\mu ^m$-almost
everywhere on $A$ an approximate differential $\mbox{ }_{ap}Df(x)$.
\par {\bf 19. Theorem.} {\it If $A\subset {\bf K}^m$, $m, n\in \bf N$,
$f: A\to {\bf K}^n$ and for each $y\in A$ there exist $\delta
=\delta (y)>0$ and $0<r=r(y)\le 1$ such that $ap {\overline
{\lim}}_{x\to z} |f(x)-f(z)|/|x-z|^r <\infty $ whenever $z\in A$ and
$|z-y|<\delta $, then $A = \bigcup_{j\in \Lambda }E_j$, where $card
(\Lambda )\le \aleph _0$, such that the restriction of $f$ to each
$E_j$ is lipschitzian; moreover, $f$ is approximately differentiable
such that there exist $\mbox{ }_{ap}{\bar {\Phi }}^1f(x;v;t)$ and
$\mbox{ }_{ap}Df(x)$ for $\mu ^{2m+1}$ almost all points $(x;v;t)$
of $A^{(1)}$ and $\mu ^m$ almost all points $x$ of $A$
correspondingly.}
\par {\bf Proof.} From the supposition of this theorem there follows
that the density of ${\bf K}^m\setminus A$ is zero and $f$ is
approximately continuous at each point of $A$. Therefore, $A$ is
$\mu ^m$ measurable and $f$ is $\mu ^m|_A$ measurable in accordance
with Theorems 2.9.11 and 2.9.13 \cite{federer}. Denote $Q_{R,j}(z):=
B({\bf K}^m,z,R)\cap \{ x: x\notin A \mbox{ or } |f(x)-f(z)|> |\pi
|^{-j}|x-z|^r \} $ for $|\pi |^j<R<\delta (y)$ and $z\in A$ with
$|y-z|<\delta =\delta (y)$ provided by the conditions of this
theorem, $0<R\in \Gamma _{\bf K}$, $j\in \bf N$. Each set $E_j :=
A\cap \{ z: \mu ^m(Q_{R,j}(z))<R^m/2 \mbox{ for } 0<R<|\pi |^j \} $
is $\mu ^m$ measurable in accordance with Lemma 10 and
$A=\bigcup_{j=1}^{\infty }E_j$. If $R=|x-z|<\delta (y)$, then $\mu
^m(Q_{R,j}(z)\cup Q_{R,j}(x))< R^m = \mu ^m(B({\bf K}^m,z,R)\cap
B({\bf K}^m,x,R))$, since $B({\bf K}^m,z,R)=B({\bf K}^m,x,R)$.
Choosing $w\in B({\bf K}^m,z,R)\setminus (Q_{R,j}(x)\cup
Q_{R,j}(z))$ we get $|f(x)-f(z)|\le \max (|f(x)-f(w)|, |f(w)-f(z)|)
\le |\pi |^j \max (|x-w|^r,|w-z|^r)\le |\pi |^j |x-z|^r$. Therefore,
if $x, z\in E_j$ and $|x-z|<|\pi |^j<R$, then $|f(x)-f(z)|\le |\pi
|^{-j} |x-z|^r$, since $A$ is everywhere dense in ${\bf K}^m$. Then
each $E_j$ is of diameter less than $|\pi |^j$, $j\in \bf N$. Each
restriction $f|_{E_j}$ is lipschitzian and Lemma 18 gives that
$f|_{E_j}$ is approximately differentiable, hence $f$ has ${\mu
}^{2m+1}$-everywhere in $A^{(1)}$ an approximate partial difference
quotient ${\bar {\Phi }}^1f(x;v;t)$ and ${\bar {\Phi }}^1f(x;e_i;t)$
for $\mu ^{m+1}$ almost all points $(x;e_i;t)$ of $(A\times \{ e_i
\} \times {\bf K})\cap A^{(1)}$ and $\mu ^m$-almost everywhere on
$A$ an approximate differential $Df(x)$ and approximate partial
differentials $D_jf(x)=\partial f(x)/\partial \mbox{ }_jx$.
\par {\bf 20. Theorem.} {\it  If $A\subset W\subset {\bf K}^m$,
$m, n \in \bf N$, $W$ is open, $A$ is $\mu ^m$ measurable, $f: W\to
{\bf K}^n$ and for each $y\in A$ there exist $\delta =\delta (y)>0$
and $0<r=r(y)\le 1$ such that ${\overline {\lim}}_{x\to z}
|f(x)-f(z)|/|x-z|^r <\infty $ whenever $z\in A$ and $|z-y|<\delta $,
then $f$ is differentiable such that there exist ${\bar {\Phi
}}^1f(x;v;t)$ and $Df(x)$ for $\mu ^{2m+1}$ almost all points
$(x;v;t)$ of $W_A^{(1)} := \{ (x;v;t): x\in A, v\in {\bf K}^m, t\in
{\bf K}, x+vt\in W \} $ and $\mu ^m$ almost all points $x$ of $A$
respectively.}
\par {\bf Proof.} The field $\bf K$ is locally compact
with the non archimedean multiplicative norm, hence $A$ has a
countable covering by balls $B({\bf K}^m,y_j,\delta _j)$ contained
in $W$, where $y_j\in A$, $\delta _j := \delta (y_j)$. Therefore,
$A$ is contained in the countable union of the subsets $E_j := W\cap
\{ z: |f(x)-f(z)| \le |\pi |^j|x-z|^r \mbox{ for }x, z\in B({\bf
K}^m,y_i,\delta _i)\mbox{ and }|\pi |^j\le \delta _i\mbox{ for some
} i\in {\bf N} \} $. Suppose that there exists a sequence $\zeta
_l\in E_j$ converging to $z\in {\bf K}^m$ as $l$ tends to the
infinity. Take $x\in B({\bf K}^m,z,|\pi |^j)\subset B({\bf
K}^m,y_i,\delta _i)$. There exists $l_0\in \bf N$ such that $\{ z,x
\} \subset B({\bf K}^m,\zeta _l,|\pi |^j)\subset W$ for each $l\ge
l_0$, hence $|f(x)-f(z)|\le \max (|f(x)-f(\zeta _l)|, |f(z)-f(\zeta
_l)|)\le |\pi |^j \max (|x - \zeta _j|^r, |z - \zeta _j|^r)\le |\pi
|^j|x-z|^r$, consequently, $z\in E_j$. Thus each $E_j$ is closed in
${\bf K}^m$. Each $E_j$ is of diameter not greater, than $\delta _i$
with the corresponding $i$ and $f|_{E_j}$ is lipschitzian. In view
of Theorems 2.9.11 \cite{federer} and 19 above the function
$f|_{E_j}$ is approximately differentiable such that there exist
$\mbox{ }_{ap}{\bar {\Phi }}^1f(x;v;t)$ and $\mbox{ }_{ap}Df(x)$ for
$\mu ^{2m+1}$ almost all points $(x;v;t)$ of $W_{E_j}^{(1)}$ and
$\mu ^m$ almost all points $x$ of $E_j$ correspondingly, since for
$t\ne 0$ we have $\mbox{ }_{ap}{\bar {\Phi
}}^1f(x;v;t)=[f(x+vt)-f(x)]/t$, $x+vt\in W$, $x\in E_j\subset A$.
Moreover, ${\bf K}^m\setminus E_j$ has zero density at each $y\in
A$. Then Theorem 8 and Lemma 16 provide that $f$ is differentiable
such that there exist ${\bar {\Phi }}^1f(x;v;t)$ and $Df(x)$ for
$\mu ^{2m+1}$ almost all points $(x;v;t)$ of $A^{(1)}$ and $\mu ^m$
almost all points $x$ of $A$ respectively.
\par {\bf 21. Lemma.} {\it  If $G\subset V\subset {\bf K}^m$,
$h: V\to {\bf K}\setminus \{ 0 \} $ is Lispchitzian with
$r=Lip_2(h)=1$, $ \{ B({\bf K}^m,y,|h(y)|): y\in G \} $ is the
disjoint family, $b\ge Lip_1(h)$, $0<\alpha \in \Gamma _{\bf K}$,
$0<\beta \in \Gamma _{\bf K}$, $b\alpha <1$ and $b\beta <1$ and
$G_x:= G\cap \{ y: B({\bf K}^m,x,\alpha |h(x)|)\cap B({\bf
K}^m,y,\beta |h(y)|) \ne \emptyset \} $ for $x\in V$, then
\par $(1)$ $[(1-b\beta)/(1+b\alpha )]\le |h(x)|/|h(y)|\le [(1+b\beta
)/(1-b\alpha )]$ for all $y\in G_x$ and \par $(2)$ $card (G_x)\le
[\max (\alpha , \beta (1+b\alpha )/(1-b\beta ))]^m[(1+b\beta
)/(1-b\alpha )]^m$.}
\par {\bf Proof.} If $y\in G_x$, then $|h(x)-h(y)|\le b|x-y|\le
b\max (\alpha |h(x)|, \beta |h(y)|)\le b(\alpha |h(x)| +\beta
|h(y)|)$, hence $(1-b\alpha )|h(x)|\le (1+b\beta )|h(y)|$ and
$(1-b\beta )|h(y)|\le (1+b\alpha )|h(x)|$, since $|x-z|_{{\bf
K}^m}\ge ||x|_{{\bf K}^m}-|z|_{{\bf K}^m}|_{\bf R}$ for each $x,
z\in {\bf K}^m$. If $y\ne g\in G$, then $B({\bf K}^m,y,|h(y)|)\cap
B({\bf K}^m,g,|h(g)|)=\emptyset $ if and only if $|y-g|>\max
(|h(y)|, |h(g)|)$. Then $|x-y|\le \max (\alpha |h(x)|, \beta
|h(y)|)\le |h(x)| \max (\alpha , \beta (1+b\alpha )/(1-b\beta ))$,
consequently, $B({\bf K}^m,y,|h(y)|)\subset B({\bf K}^m,x,\gamma
|h(x)|)$, where $\gamma := \max (\alpha , \beta (1+b\alpha
)/(1-b\beta ))$, hence $card (G_x) [(1-b\alpha )|h(x)|/(1+b\beta
)]^m\le \sum_{y\in G_x} |h(y)|^m\le (\beta (1+b\alpha )/(1-b\beta
))^m|h(x)|^m$, since $G_x\subset G$ and $G$ is discrete, $G_x\subset
B({\bf K}^m,x,\gamma |h(x)|)$.
\par {\bf 22. Theorem.} {\it  Suppose $Y$ is a normed vector space over
$\bf K$, $A$ is a closed subset in ${\bf K}^m$ and to each $z\in A$
there corresponds a polynomial function $P_z: {\bf K}^m\to Y$ with
degree $deg (P_z)\le k$. Let also $S\subset A$ and $\delta >0$ and
$\rho (S,\delta ):=\sup_{0<|x-z|\le \delta ; x, z\in S; j=0,1,...,k}
\| {\bar {\Phi }}^jP_x(z;v;t) - {\bar {\Phi }}^jP_z(z;v;t)
\|_{C^0(V_z^{(j)},Y)} |x-z|^{j-k}$, where $U_z=B({\bf K}^m,z,|\pi
|^{\zeta })$, $\zeta \in \bf N$ is a marked number,
$v=(v_1,...,v_j)$, $t=(t_1,...,t_j)$, $U^{(j)}:=\{
(y;v_1,...,v_j;t_1,...,t_j): y\in U, y+v_1t_1+...+v_jt_j\in U,
v_i\in {\bf K}^m, t_i\in {\bf K} \quad \forall i=1,...,j \} $,
$V_z^{(j)}:= \{ (y;v_1,...,v_j;t_1,...,t_j)\in U_z^{(j)}: v_i\in
{\bf K}^m, |v_i|=1, t_i\in {\bf K} \quad \forall i=1,...,j \} $. If
$\lim_{0<\delta \to 0} \rho (S,\delta )=0$ for each compact subset
$S$ of ${\bf K}^m$, then there exists a map $g: {\bf K}^m\to Y$ of
class $C^k$ such that ${\bar {\Phi }}^jg(z;v;t)={\bar {\Phi
}}^jP_z(z;v;t)$ on $V_z^{(j)}$ for each $j=0,1,...,k$ and $z\in A$.}
\par {\bf Proof.} Let ${\cal F}$ be a family of open subsets of
${\bf K}^m$. Take $b = |\pi |^{s_0}$ for a marked natural number
$s_0$. Let $h_R$ be a function on $W:=\bigcup_{U\in \cal F} U$ such
that $h_R(x) := b \sup \{ \inf (1, dist (x,{\bf K}^m\setminus T)):
T\in {\cal F} \} $. Since $\Gamma _{\bf K}$ is discrete in
$(0,\infty )$ and $|*|$ is the continuous norm from $\bf K$ into
$\Gamma _{\bf K}\cup \{ 0 \} $, then $h_R(x)\in \Gamma _{\bf K}\cup
\{ 0 \} $ for each $x\in W$ and $h_R(x)$ is continuous, hence for
each clopen or closed $G$ in $\Gamma _{\bf K}\cup \{ 0 \} $ its
counter image $h_R^{-1}(G)$ is clopen or closed in $W$ respectively.
Therefore, $W$ is the disjoint union of the closed set $h_R^{-1}(0)$
and the clopen subsets $h_R^{-1}(u)$ while $u\in \Gamma _{\bf K}$,
consequently, there exists a continuous function $h: W\to \bf K$
such that $|h(x)|=h_R(x)$ for each $x\in W$, since $h_R$ is
continuous.
\par In accordance with Theorem 2.8.4 \cite{federer}
if $(X,\rho _X)$ is a metric space and $F$ is a family of its closed
subsets, $\delta $ is a nonnegative bounded function on $F$ and
$1<\tau <\infty $, then $F$ has a disjointed subfamily $G$ such that
for each $T\in F$ there exists $H\in G$ with $T\cap H\ne \emptyset $
and $\delta (T)\le \tau \delta (H)$.
\par With each $H\in F$ it is possible to associate its $\delta ,
\tau $ enlargement ${\hat H}:=\bigcup \{ T: T\in F, T\cap H\ne
\emptyset , \delta (T)\le \tau \delta (H) \} $. Corollary 2.8.5
\cite{federer} states that $\bigcup_{T\in F} T\subset \bigcup_{H\in
G}{\hat H}$.
\par In the considered here situation take $F:= \{ B({\bf
K}^m,x,h_R(x)): x\in W \} $, $\delta =diam $, $\tau = |\pi |^{s_1}$
for a marked integer number $s_1\le 0$. Choose $G\subset W$ so that
$ \{ B({\bf K}^m,y,h_R(y)): y\in  G \} $ is disjoined, $h_R(y)>0$,
and $\bigcup_{y\in  G} B({\bf K}^m,y,h_R(y))=W$, where $s_0, s_1$
are subordinated to the condition $|s_1|+1<s_0$. Evidently $G$ is
countable. By Lemma 21 we get that $x\in W$ implies $(1-b\beta
)/(1+b\alpha )\le h_R(x)/h_R(y)\le (1+b\beta )/(1-b\alpha )$ for all
$y\in G_x$. Taking $\alpha =\beta =b|\pi |^{s_2}$ we get $[(1-|\pi
|^{2s_0+s_2})/(1+|\pi |^{2s_0+s_2}]< |h(x)|/|h(y)|\le [(1+|\pi
|^{2s_0+s_2})/(1-|\pi |^{2s_0+s_2}]$ for every $y\in G_x$. Moreover,
$card (G_x)\le |\pi |^{s_0+s_2}[(1+|\pi |^{2s_0+s_2})/(1-|\pi
|^{2s_0+s_2})]^{2m}$, where $s_0\ge 1$, $s_2\ge -1$ are integers.
\par Consider the mapping $w_y(x):=ch_{B({\bf K}^m,0,1)}((x-y)/(\pi
h(y))$ for $y\in G$, $x\in {\bf K}^m$, where $ch_P$ is the
characteristic function of a subset $P$ in ${\bf K}^m$. Therefore,
$supp w_y = B({\bf K}^m,y,|\pi h(y)|)$. This function is of
$C^{\infty }$ class with finite $C^n$ norms for each $n\in \bf N$,
since each ${\bar {\Phi }}^jw_y(x;v;t)$ is bounded on ${\bf
K}^m\times S({\bf K}^m,0,1)^j\times {\bf K}^j$ for each $j\in \bf
N$, where $S({\bf K}^m,z,R) := \{ x\in {\bf K}^m: |x-z|=R \} $ for
$R\in \Gamma _{\bf K}$, $|x|= \max_{j=1}^m|\mbox{ }_jx|$, $x=(\mbox{
}_1x,...,\mbox{ }_mx)$. Indeed, ${\bar {\Phi }}^1ch_B(x;v;t)=0$ for
$x, x+vt\in B$ or for $|x|>1$ and $|x+vt|>1$, ${\bar {\Phi
}}^1ch_B(x;v;t)=1/t$ for $|x|\le 1$ and $|x+vt|>1$ or $|x|>1$ and
$|x+vt|\le 1$, where $B=B({\bf K}^m,0,1)$. In the latter two cases
$|t|>1$ for $|v|=1$. Continuing by induction we get $|{\bar {\Phi
}}^jch_{B({\bf K}^m,0,1)}(x;v;t)|\le 1$ for each $x\in {\bf K}^m$,
$v\in S({\bf K}^m,0,1)^j$ and $t \in {\bf K}^j$, hence $|{\bar {\Phi
}}^jw_y(x;v;t)|\le |\pi h(y)|^{-j}$ for each $(x;v;t)\in {\bf
K}^m\times S({\bf K}^m,0,1)^j\times {\bf K}^j$. \par Choose now
$G_0\subset G$ with the additional condition that \par
$\bigcup_{y\in G_0}B({\bf K}^m,y,|\pi h(y)|)\supset W$ and \par
$B({\bf K}^m,y,|\pi h(y)|)\cap B({\bf K}^m,g,|\pi h(g)|)=\emptyset $
for each $y\ne g\in G_0$\\ and consider the function $\phi
(x):=\sum_{y\in G_0}w_y(x)$. Then $\phi \in C^{\infty }$ and $\phi
(x)|_W=1$ for each $x\in W$. Thus the family of functions $\{
w_y(x): y\in G_0 \} $ constitute the partition of unity on $W$
associated with the family $\cal F$. They are of class $C^{\infty }$
and their supports form a disjoint clopen refinement $\cal F$ of
covering $W$.
\par Consider the subset $U={\bf K}^m\setminus A$ and put ${\cal F}=
\{ U \} $. Since $A$ is closed, then $U$ is open in ${\bf K}^m$.
Applying Lemma 21 we get $h_R(x)/b=\inf \{ 1, dist (x,A) \} $ for
$x\in U$. For $y\in S$ we take $\psi (y)\in A$ with $|y-\psi (y)|=
dist (y,A)$. Then we define a function $g: {\bf K}^m\to Y$ by the
formula $g(x)=P_x(x)$ for $x\in A$ and $g(x)=\sum_{y\in S}
w_y(x)P_{\psi (x)}(x)$ for $x\in U$. Therefore, $g\in C^{\infty
}({\bf K}^m,{\bf K})$, since ${\bar {\Phi }}^jg$ are polynomials of
${\bar {\Phi }}^iw_y$ and ${\bar {\Phi }}^lP_{\psi (y)}$ with the
corresponding arguments, where $1\le i\le j$, $1\le l\le j$ (see
Corollary 2.6). \par In the particular case of $X={\bf K}^m$
we have \par ${\bar {\Phi }}^jf(z;y-z,...,y-z;0,...,0)=$ \\
$\sum_{l(1),...,l(j)\in \{ 1,...,m \} } {\bar
{\Phi}}^jf(z;e_{l(1)},...,e_{l(j)};0,...,0) (\mbox{ }_{l(1)}y-\mbox{
}_{l(1)}z)... (\mbox{ }_{l(j)}y-\mbox{ }_{l(j)}z)$
\\ for $f\in C^n(U,Y)$, $j\le n$, where $U$ is an open subset in
${\bf K}^m$, $y, z\in U$. In view of Theorem A.1 (see Appendix)
\par ${\bar {\Phi }}^jP_x(y^{(j)}) - {\bar {\Phi }}^jP_z(y^{(j)}) =
\sum_{i=0}^k\sum_{{\bar l}_i} ({\bar {\Phi }}^i_y[{\bar {\Phi
}}^jP_x(y^{(j)}) - {\bar {\Phi }}^jP_z(y^{(j)})](z^{(i+j)}_{{\bar
l}_i})).(y-z)^i + \sum_{{\bar l}_{k-j}}R_{k-j}({\bar {\Phi
}}^jP_x(z^{(j)}) - {\bar {\Phi }}^jP_z(z^{(j)});z^{(k)}_{{\bar
l}_{k-j}}).(y-z)^{k-j}$, \\ where $x, z\in S$, $y\in {\bf K}^m$,
$j\le k$, $y^{(j)}=(y;v_1,...,v_j;t_1,...,t_j)$, $z^{(i+j)}_{{\bar
l}_i}=(z;v_1,...,v_j,e_{l(1)},...,e_{l(i)};t_1,...,t_j,0,...,0)$,
$v_i\in {\bf K}^m$, $t_i\in \bf K$, ${\bar
l}_i=(l(1),...,l(i))\subset \{ 1,...,m \} $ for each $i=1,...,j$;
\par $({\bar {\Phi }}^i_y[{\bar {\Phi
}}^jP_x(y^{(j)})](z^{(i+j)}_{{\bar l}_i})).(y-z)^i$ \\  $= {\bar
{\Phi }}^{i+j}P_x(z^{(i+j)}_{{\bar l}_i}) (\mbox{ }_{l(1)}y-\mbox{
}_{l(1)}z)... (\mbox{ }_{l(i)}y-\mbox{ }_{l(i)}z)$; \\
$R_{k-j}({\bar {\Phi }}^jP_x(z^{(j)}) - {\bar {\Phi
}}^jP_z(z^{(j)});z^{(k)}_{{\bar l}_{k-j}})$  is the continuous
residue equal to zero for $x=z$ for each ${\bar l}_{k-j}$,
\par $R_{k-j}(f;z^{(j)}_{{\bar l}_{k-j}}).(y-z)^{k-j}:=
R_{k-j}(f;z^{(j)}_{{\bar l}_{k-j}})(\mbox{ }_{l(1)}y-\mbox{
}_{l(1)}z)... (\mbox{ }_{l(k-j)}y-\mbox{ }_{l(k-j)}z)$.
\par Therefore,
\par $\| {\bar {\Phi }}^jP_x(y^{(j)}) - {\bar {\Phi }}^jP_z(y^{(j)})\|
\le \max _{0\le i\le k-j} (|y-z|^i|x-z|^{k-j-i}) \rho (S,|x-z|)$.
\par If $z\in A$, $G=A\cap B({\bf K}^m,z,|\pi |^{-1})$ and
$x\in U\cap B({\bf K}^m,z,|\pi |)$, then we choose $y\in G$ with
$|x-y|=dist (x,A)$, hence $|x-y|\le |x-z|\le |\pi |^{-1}$, $|\pi
|^{s_0}|h(x)|=|x-y|\le |\pi |$ and $|y-z|\le \max (|x-y|, |x-z|)\le
|x-z|<|\pi |^{-1}$. Take $s_2\ge -1$ such that $[(1+|\pi
|^{2s_0+s_2})/(1-|\pi |^{2s_0+s_2})]\le |\pi |^{-1}$, consequently,
$q\in G_x$ implies $|\pi |^{-s_0}|h(q)|\le |\pi |^{-s_0}[(1+|\pi
|^{2s_0+s_2})/(1-|\pi |^{2s_0+s_2})]|h(x)|\le 1$ and $|\pi |^{-s_0}
|h(q)|=|q-\psi (q)|$, $|q-x|\le \max (|h(q)|,|h(x)|)|\pi |^{s_0+s_2}
\le |\pi |^{s_0+s_2}\le 1$ and $|\psi (q)-z|\le \max (|\psi (q)-q|,
|q-x|,|x-z|)\le \max (|\pi |^{-1},1,|\pi |^{-1})=|\pi |^{-1}$ and
$|\psi (q)-y|\le \max (|\psi (q)-q|,|q-x|,|x-y|)\le \max (|\pi
|^{-1},1,|\pi |^{-1})=|\pi |^{-1}$, since $s_0+s_2\ge 0$, where
$\psi (q)\in G$. \par If $f\in C^n(U,{\bf K})$ and $g\in C^n(U,Y)$,
where $U$ is open in ${\bf K}^m$, then
\par $(i)$ ${\bar {\Phi }}^n(fg)(x^{(n)}) =
\sum_{0\le a, 0\le b, a+b=n}\sum_{j_1<...<j_a; s_1<...<s_b; \{
j_1,...,j_a \} \cup \{ s_1,...,s_b \} = \{ 1,...,n \} }$ \\
${\bar {\Phi }}^af(x;v_{j_1},...,v_{j_a};t_{j_1},...,t_{j_a}) {\bar
{\Phi }}^bg(x+v_{j_1}t_{j_1}+...+v_{j_a}t_{j_a};v_{s_1},...,v_{s_b};
t_{s_1},...,t_{s_b})$ (see also Corollary 2.6).
\par Take $U_z=U_y=B({\bf K}^m,z,|\pi |^{\zeta })$, that is, $y\in B({\bf
K}^m,z,|\pi |^{\zeta })$ with $\zeta \in \bf N$, then
\par $ \| {\bar {\Phi }}^ig(x^{(i)})-{\bar {\Phi }}^iP_y(x^{(i)}) \|
_{C^0(V^{(i)}_y,{\bf K})}  \le \sup_{q\in G_x}\max_{0\le j\le i} \|
{\bar {\Phi }}^{i-j}w_q \| \max_{0\le l\le j}(|x-y|^i|\psi
(q)-y|^{j-i})\rho (G,|\psi (q)-y|)\le |\pi |^{-i}\rho (G,|\psi
(q)-y|)$, since ${\bar {\Phi }}^ig(x^{(i)})-{\bar {\Phi
}}^iP_y(x^{(i)}) = \sum_{q\in G} {\bar {\Phi }}^i(w_q(P_{\psi
(q)}-P_y))(x^{(i)})$ and due to Formula $(i)$, also $\max_{0\le j\le
i} \| {\bar {\Phi }}^{i-j}w_q \| \le \max_{0\le j\le i} |\pi
h(q)|^{i-j}\le \max_{0\le j\le i} |\pi |^{(s_0+1)(i-j)}\le 1$. Thus,
\par $(ii)$ $ \| {\bar {\Phi }}^ig(x^{(i)})-{\bar {\Phi }}^iP_z(x^{(i)}) \|
_{C^0(V^{(i)}_z,{\bf K})} \le $\\ $\max ( \| {\bar {\Phi
}}^ig(x^{(i)})-{\bar {\Phi }}^iP_y(x^{(i)}) \| _{C^0(V^{(i)}_z,{\bf
K})}, \| {\bar {\Phi }}^iP_y(x^{(i)})-{\bar {\Phi }}^iP_z(x^{(i)})
\| _{C^0(V^{(i)}_z,{\bf K})}\le |\pi |^{-i} \max (\rho (G,|\psi
(q)-y|), \rho (G,|y-z|))$ \\
for each $i\le k$. Therefore, by induction relative to $i$ we get
${\bar {\Phi }}^ig(x^{(i)})={\bar {\Phi }}^iP_x(x^{(i)})$ for each
$x\in A$. From Formula $(ii)$ it follows that there exists
\par $\lim_{x\to z} \| {\bar {\Phi }}^ig(x^{(i)})-{\bar {\Phi
}}^iP_z(x^{(i)}) \|_{C^0(V^{(i)}_z,{\bf K})}|x-z|^{i-k}=0$ \\
for all $z\in A$ and $1\le i\le k$, hence ${\bar {\Phi }}^ig$ is
continuous at $z$ and for $i<k$ and $1\le l\le k-i$ there are the
identities ${\bar {\Phi }}^l({\bar {\Phi
}}^ig(z^{(i)}))(z^{i+l})={\bar {\Phi }}^l({\bar {\Phi
}}^iP_z(z^{(i)}))(z^{i+l})={\bar {\Phi }}^{i+l}P_z(z^{(i+l)})={\bar
{\Phi }}^{i+l}g(z^{(i+l)})$, where ${\bar {\Phi }}^ig(z^{(i)})$ and
${\bar {\Phi }}^iP_z(z^{(i)})$ are defined on $U^{(i)}_z$ such that
$z\in U_z$, $U_z$ is open in ${\bf K}^m$ and it is sufficient to
consider $U_z=B({\bf K}^m,z,|\pi |^{\zeta })$.
\par {\bf 23. Theorem.} {\it If $A\subset W\subset {\bf K}^m$,
$f\in C^k(W,{\bf K}^n)$ and \par ${\overline {\lim}}_{x^{(k)}\to
z^{(k)}} \| {\bar {\Phi }}^kf(x^{(k)})-{\bar {\Phi }}^kf(z^{(k)}) \|
_{C^0(V^{(k)}_{z,R},{\bf K}^n)}/|x^{(k)}-z^{(k)}|^r<\infty $ \\
for each $z\in A$ and each $0<R<\infty $, where $W$ is open in ${\bf
K}^m$ and $0<r\le 1$, $V^{(k)}_{z,R}$ corresponds to $U_{z,R}=W\cap
B({\bf K}^m,z,R)$, then for each $\epsilon
>0$ there exists a map $g\in C^{k+1}({\bf K}^m,{\bf K}^n)$ such that
$\mu ^m(A\setminus \{ x: f(x)=g(x) \} )<\epsilon $.}
\par {\bf Proof.} In accordance with Theorem 20 $\mu
^{m+(m+1)(k+1)}(A^{(k+1)}\setminus dom ({\bar {\Phi }}^{k+1}f))=0$
and due to Theorem 15 ${\bar {\Phi }}^{k+1}f$ is $\mu
^{m+(m+1)(k+1)}|_{dom ({\bar {\Phi }}^{k+1}f)}$ measurable. From the
Lousin Theorem 2.3.5 \cite{federer} it follows, that there exists a
closed subset $E$ in $A^{(k+1)}$ such that ${\bar {\Phi
}}^{k+1}f(x^{(k+1)})|_E$ is continuous and $\mu
^{m+(m+1)(k+1)}(A^{(k+1)}\setminus E)<\epsilon $, where $A^{(k+1)}$
is defined analogously to $U^{(k+1)}$. Practically ${\bar {\Phi
}}^{k+1}f(x^{(k+1)})$ is also continuous on $W^{(k)}\times {\bf K}^m
\times {\bf K}\setminus \{ 0 \} $, that is, for $t_{k+1}\ne 0$. Put
\par $\phi _q(z^{(k)}) :=
\sup_{0\le j\le k} \{ \| {\bar {\Phi }}^jf(x^{(j)})- {\bar {\Phi
}}^jf(z^{(j)})- {\bar {\Phi }}^1({\bar {\Phi
}}^jf(z^{(j)}))(z^{(j)};x^{(j)}-z^{(j)};1) \| /|x^{(j)}-z^{(j)}|^r:
x^{(j)}\in V^{(j)}_{z,|\pi
|^q}; x^{(j)}\ne z^{(j)} \} $, \\
where $z^{(k)}=(z;v_1,...,v_k;t_1,...,t_k)\in A^{(k)}$,
$z^{(j)}=(z;v_1,...,v_j;t_1,...,t_j)$, $q=1,2,...$. Then there
exists $\lim_{q\to \infty } \phi _q(z^{(k)})=0$ for each $z\in E$.
Each function $\phi _q$ is borelian.
\par Theorem 2.2.2 \cite{federer} states: suppose nonnegative
$\phi $ is a measure over a metric space $X$, all open subsets of
$X$ are $\phi $ measurable, and $B$ is a Borel set; $(1)$ if $\phi
(B)<\infty $ and $\epsilon >0$, then $B$ contains a closed set $C$
for which $\phi (B\setminus C)<\epsilon $; $(2)$ if $B$ is contained
in the union of countably many open sets $V_i$ with $\phi
(V_i)<\infty $, and if $\epsilon >0$, then $B$ is contained in an
open set $W$ for which $\phi (W\setminus B)<\epsilon $.
\par While the Egoroff Theorem 2.3.7 \cite{federer} is:
suppose $f_1, f_2,...$ and $g$ are $\phi $ measurable functions with
values in a separable metric space $Y$, where nonnegative $\phi $ is
a measure over $X$; if $\phi (A)<\infty $, $A\subset X$, $f_n(x)\to
g(x)$ for almost all $x$ in $A$, and $\epsilon >0$, then there
exists a $\phi $ measurable set $B$ such that $\phi (A\setminus
B)<\epsilon $ and $f_n(x)\to g(x)$, uniformly for $x\in B$, as $n\to
\infty $.
\par In view of the Egoroff 2.3.7 and 2.2.2 Theorems \cite{federer} for
each $\epsilon
>0$ there exists a closed subset $F$ in $A^{(k)}$ such that $\mu
^{m+(m+1)k)}(A^{(k)}\setminus F)<\epsilon $ and $\lim_{q\to \infty }
\sup_{z\in J}\phi _q(z)=0$ for each compact subset $J$ in $F$.
\par Take $G$ and $G_{\epsilon }$ from Theorem 17 and for $z\in G_{\epsilon }$
consider the Taylor expansion of Theorem A.1 with $k$ here instead
of $n$ there and put $P_z(y)=f(z)+\sum_{j=1}^{k+1} {\bar {\Phi
}}^jf(z;y-z,...,y-z;0,...,0)$, hence $P_y(y)=f(y)$. On the other
hand $\mu ^m(A\setminus G)=0$ due to Theorems 17 and 19. Show that
the suppositions of Theorem 22 are satisfied with $Y={\bf K}^n$. Let
$H$ be a compact subset in $G_{\epsilon }$ and $y, z\in H$. In
accordance with Theorem A.1 we have:
\par ${\bar {\Phi }}^iP_y(y^{(i)}) - {\bar {\Phi }}^iP_z(y^{(i)})=
{\bar {\Phi }}^if(y^{(i)})- {\bar {\Phi }}^if(z^{(i)})-
[\sum_{j=i+1}^{k+1}{\bar {\Phi
}}^jf(z;v_1,...,v_i,y-z,...,y-z;t_1,...,t_i,0,...,0)] $ \\
for each $i=0,...,k-1$ and $0<|y-z|<|\pi |^q$, where
$y^{(i)}=(y;v_1,...,v_i;t_1,...,t_i)$ and
$z^{(i)}=(z;v_1,...,v_i;t_1,...,t_i)$. Therefore, \\
$\sup_{0<|y-z|\le |\pi |^q ; y, z\in H; 0\le i\le k} | {\bar {\Phi
}}^iP_y(y^{(i)})- {\bar {\Phi }}^iP_z(y^{(i)})|
|y-z|^{i-k-1}\le \phi _q(z^{(k)})<\infty $, \\
where $z^{(k)}=(z;v_1,...,v_i,y-z,...,y-z;t_1,...,t_i,0,...,0)$.
Since ${\bar {\Phi }}^{k+1}f$ is continuous on $G$, then ${\bar
{\Phi }}^{k+1}P_y(y^{(k+1)})-{\bar {\Phi }}^{k+1}P_z(y^{(i)})= {\bar
{\Phi }}^{k+1}f(y^{(k+1)})-{\bar {\Phi }}^{k+1}f(z^{(k+1)})$ is
small for small $|y^{(k+1)}-z^{(k+1)}|$. Thus from Theorem 22 the
statement of this theorem follows.
\par {\bf 24. Theorem.} {\it 1. If $A\subset {\bf K}^m$, $f: A\to {\bf
K}^n$ and $ap {\overline {\lim}}_{x\to z} |f(x)-f(z)|/|x-z|^r$ for
$\mu ^m$ almost all $z\in A$, where $0<r\le 1$, then for each
$\epsilon >0$ there exists a map $g\in C^1({\bf K}^m,{\bf K}^n)$
such that $\mu ^m(A\setminus \{ x: f(x)=g(x) \} )<\epsilon $. \par
2. If $A\subset W\subset {\bf K}^m$, $f\in C^k(W,{\bf K}^n)$ and
\par $ap {\overline {\lim}}_{x^{(k)}\to z^{(k)}} \| {\bar {\Phi
}}^kf(x^{(k)})-{\bar {\Phi }}^kf(z^{(k)}) \|
_{C^0(V^{(k)}_{z,R},{\bf K}^n)}/|x^{(k)}-z^{(k)}|^r<\infty $ \\
for $\mu ^m$ almost all $z\in A$ and each $0<R<\infty $, where $W$
is open in ${\bf K}^m$ and $0<r\le 1$, $A$ is $\mu ^m$ measurable,
$V^{(k)}_{z,R}$ corresponds to $U_{z,R}=W\cap B({\bf K}^m,z,R)$,
$k\ge 1$, then for each $\epsilon >0$ there exists a map $g\in
C^{k+1}({\bf K}^m,{\bf K}^n)$ such that $\mu
^{m+(m+1)k}(W_A^{(k)}\setminus \{ x^{(k)}: {\bar {\Phi
}}^kf(x^{(k)})={\bar {\Phi }}^kg(x^{(k)}) \} )<\epsilon $, where
$W_A^{(k)} := \{ (x;v_1,...,v_k;t_1,...,t_k): x\in A; v_1,...,v_k\in
{\bf K}^m; t_1,...,t_k\in {\bf K}; x+v_1t_1+...+v_jt_j\in W \quad
\forall 1\le j\le k \} $.}
\par {\bf Proof.} In accordance with Theorems 19 and 20
there exist compact subsets $H_1, H_2,... $ such that $H_i\subset
B({\bf K}^{m+(m+1)k},0,|\pi |^{-i})\setminus B({\bf
K}^{m+(m+1)k},0,|\pi |^{-i+1})$ for $i>1$ and $H_1\subset B({\bf
K}^{m+(m+1)k},0,1/|\pi |)$ and the restrictions ${\bar {\Phi
}}^kf|_{H_i}$ are lipschitzian, $\mu ^{m+(m+1)k}(A^{(k)}\cap
[(B({\bf K}^{m+(m+1)k},0,|\pi |^{-i})\setminus B({\bf
K}^{m+(m+1)k},0,|\pi |^{-i+1}))\setminus H_i])<\epsilon 2^{-i}$ for
$i>1$ and $\mu ^{m+(m+1)k}(W_A^{(k)}\cap [B({\bf
K}^{m+(m+1)k},0,1/|\pi |)\setminus H_1])<\epsilon /2$, where
$W_A^{(0)}=A$, ${\bar {\Phi }}^0f=f$, $k=0$ in the first case and
$k>0$ in the second case. Using Theorem 8 by induction we construct
a function $h: {\bf K}^{m+(m+1)k}\to {\bf K}^n$ such that
$h|_{B({\bf K}^{m+(m+1)k},0,|\pi |^{-i})}$ is a lipschitzian
extension of $h|_{B({\bf K}^{m+(m+1)k},0,|\pi |^{-i+1})}\cup {\bar
{\Phi }}^kf|_{H_i}$. Therefore, $h$ is locally lipschitzian and $\mu
^{m+(m+1)k}(W_A^{(k)}\setminus \{ x^{(k)}: {\bar {\Phi
}}^kf(x^{(k)})=h(x^{(k)}) \} )<\epsilon $. The applying Theorem 23
with ${\bf K}^m, {\bf K}^n, h, k$ instead of $A, W, f, k$ gives the
assertion of this theorem.
\par {\bf 25. Theorem.} {\it Let $f: {\bf K}^m\to {\bf K}^n$,
$m, n\in \bf N$. Let also $f\circ u\in C^{s,r}({\bf K},{\bf K}^n)$
for each $u\in C^{\infty }({\bf K},{\bf K}^m)$, where $s$ is a
nonnegative integer, $0< r\le 1$, then there exists a $\mu ^m$
measurable subset $G$ in ${\bf K}^m$ such that $f\in C^{s+1}(G,{\bf
K}^n)$, where $\mu ^m({\bf K}^m\setminus G)=0$. Moreover, for each
$\epsilon >0$ there exists a map $g\in C^{s+1}({\bf K}^m,{\bf K}^n)$
such that $\mu ^m({\bf K}^m\setminus \{ x: f(x)=g(x) \} )<\epsilon
$.}
\par {\bf Proof.} The non-archimedean modification of the results from
\cite{boman} gives, that $f\in C^{s,r}({\bf K}^m,{\bf K}^n)$. Then
by Theorem 9 there exists ${\bar {\Phi }}^{s+1}f(x^{(s+1)})$
continuous for almost all points $x^{(s+1)}\in {\bf K}^m\times ({\bf
K}^m)^{s+1}\times ({\bf K})^{s+1}$ on a set $G_{0,s+1}$ and
$D^i{\bar {\Phi }}^jf(x^{(j)})$ is continuous for almost all points
of ${\bf K}^m\times ({\bf K}^m)^j\times ({\bf K})^j$ on a set
$G_{i,j}$ for each $0\le i, j\le s$ with $i+j\le s+1$. Since
$G_{i,j}^{(i)}\subset G_{0,i+j}$ for each $i, j$ and $G_{0,k}={\bf
K}^m\times ({\bf K}^m)^k\times ({\bf K})^k$ for each $k\le s$, then
$G\supset G_{0,s+1}$ and $\mu ^m({\bf K}^m\setminus G)=0$.
\par The second statement follows from Theorem 23.
\par {\bf 26. Theorem.} {\it 1. Suppose $f: {\bf K}^m\to {\bf K}^n$,
$m, n\in \bf N$ and $f\circ u\in \mbox{ }_{ap}C^{s+1}({\bf K},{\bf
K}^n)$ for each $u\in C^{\infty }({\bf K},{\bf K}^m)$, where $s$ is
a nonnegative integer, then $f\in \mbox{ }_{ap}C^{s+1}({\bf
K}^m,{\bf K}^n)$. \par 2. If $A\subset {\bf K}$, $f\circ u\in \mbox{
}_{ap}C^{s,r}({\bf K},A,{\bf K}^n)$ for each $u\in C^{\infty }({\bf
K},{\bf K}^m)$, where $0<r\le 1$, $A$ is $\mu $ measurable, then for
each $\epsilon >0$ there exists a map $g\in C^{s+1}({\bf K}^m,{\bf
K}^n)$ such that $\mu ^{m+(m+1)s)}(({\bf K}^m_{A^m})^{(s)}\setminus
\{ x^{(s)}: {\bar {\Phi }}^sf(x^{(s)})={\bar {\Phi }}^sg(x^{(s)}) \}
)<\epsilon $.}
\par {\bf Proof.} From the ultra-metric modification of
\cite{boman} it follows, that $f\in C^s({\bf K}^m,{\bf K}^n)$. \par
1. It remains to prove, that ${\bar {\Phi }}^sf \in \mbox{
}_{ap}C^1(({\bf K}^m)^{(s)},{\bf K}^n)$, where $({\bf
K}^m)^{(s)}={\bf K}^m\times ({\bf K}^m)^s\times ({\bf K})^s$ and
${\bar {\Phi }}^sf\circ u \in \mbox{ }_{ap}C^1(({\bf K})^{(s)},{\bf
K}^n)$ for each $u\in C^{\infty }({\bf K},{\bf K}^m)$. Therefore, it
is sufficient to prove, assertion 1 for $s=0$ up to a choice of
notation. But $f\circ u\in \mbox{ }_{ap}C^1({\bf K},{\bf K}^n)$ for
each $u\in C^{\infty }({\bf K},{\bf K}^m)$ means in particular this
for $u(t_0)=x$, $u(t_1)-u(t_0)=e_it$ with $t_0\ne t_1\in \bf K$,
that ${\bar {\Phi }}^1 f(x;e_i;t)$ is $\mu ^{m+1}$ almost everywhere
continuous on ${\bf K}^m\times \bf K$ and there exists a linear
mapping $T_i: {\bf K}\to {\bf K}^n$ such that $(\mu ,V) ap
\lim_{t\to 0} |f(x+e_it)-f(x)- T_it |/ |t|=0$ for each $i=1,...,m$,
where $V = {\cal B}({\bf K})$. Since $i$ is arbitrary, then ${\bar
{\Phi }}^1 f(x;v;t)$ is $\mu ^{2m+1}$ almost everywhere on ${\bf
K}^{2m+1}$ continuous due to the combinatorial Formula:
\par  ${\bar {\Phi }}^1(f\circ u)(y,v,t) = \sum_{j=1}^m
{\hat S}_{j+1,vt} {\bar {\Phi }}^1f(u(y),e_j,t{\bar {\Phi }}^1\circ
p_ju(y,v,t)) ({\bar {\Phi }}^1\circ p_ju(y,v,t))$, \\
where  $S_{j,\tau }u(y):=(u_1(y),...,u_{j-1}(y),u_j(y+\tau
_{(s)}),u_{j+1}(y+ \tau _{(s)}),...,u_m(y+\tau _{(s)}))$,
$u=(u_1,...,u_m)$, $u_j\in \bf K$ for each $j=1,...,m$, $y\in {\bf
K}^s$, $\tau =(\tau _1,...,\tau _k)\in {\bf K}^k$, $k\ge s$, $\tau
_{(s)}:=(\tau _1,...,\tau _s)$, $p_j(x):=x_j$, $x=(x_1,...,x_m)$,
$x_j\in \bf K$ for each $j=1,...,m$, ${\hat S}_{j+1,\tau
}g(u(y),\beta ):=g(S_{j+1,\tau }u(y),\beta )$, $y\in {\bf K}^s$,
$\beta $ is some parameter. Then there exists a linear mapping $T:
{\bf K}^m\to {\bf K}^n$ such that $(\mu ^m,V) ap \lim_{z\to x}
|f(z)-f(x)- T(z-x)|/ |z-x|=0$, where $V = {\cal B}({\bf K}^m)$,
since $\mu ^m=\otimes_{i=1}^m\mu $ and ${\cal B}({\bf
K}^m)=\otimes_{i=1}^m{\cal B}({\bf K})$ is the minimal $\sigma $
algebra generated from the family of all subsets of the form
$A_1\times ... \times A_m$ with $A_i\in {\cal B}({\bf K})$ for each
$i=1,...,m$, consequently, $f\in \mbox{ }_{ap}C^1({\bf K}^m,{\bf
K}^n)$, where $Te_it=T_it$ for each $t\in \bf K$ and each
$e_i=(0,...,0,1,0,...,0)\in {\bf K}^m$.
\par 2. We will prove, that ${\bar {\Phi
}}^sf \in \mbox{ }_{ap}C^{0,r}(({\bf K}^m)^{(s)},A^m,{\bf K}^n)$,
then the second statement will follow from Theorem 24.2. Thus up to
the notation it is sufficient to prove the second assertion for
$s=0$. We have $f\circ u\in \mbox{ }_{ap}C^{0,r}({\bf K},A,{\bf
K}^n)$ for each $u\in C^{\infty }({\bf K},{\bf K}^m)$, particularly,
for $u(t_0)=z$, $u(t_1)-u(t_0)=e_it$ with $t_0\ne t_1\in \bf K$.
Thus
\par $(\mu ,V) ap {\overline {\lim}}_{t\to 0} \| f(z+e_it)-f(z) \|
_{C^0(W_{t_0,R},{\bf K}^n)}/|t|^r<\infty $ \\
for $\mu $ almost all $t_0\in A$ and each $0<R<\infty $, where
$W_{t_0,R}=B({\bf K},t_0,R)$, $V={\cal B}({\bf K})$, consequently,
\par $(1)$ $(\mu ^m,V) ap {\overline {\lim}}_{x\to z} \| f(z)-f(x) \|
_{C^0(U_{z,R},{\bf K}^n)}/|x-z|^r\le \max_{i=1}^m ap {\overline
{\lim}}_{\max_i|t_i|\to 0} \|
f(z+e_1t_1+...+e_it_i)-f(z+e_1t_1+...+e_{i-1}t_{i-1}) \|
_{C^0(W_{t_{0,i},R},{\bf K}^n)}/|t_i|^r <\infty $ \\
for $\mu ^m$ almost all $t_0=(t_{0,1},...,t_{0,m})\in A^m$ and each
$0<R<\infty $, where $V={\cal B}({\bf K}^m)$, since $\mu
^m=\otimes_{i=1}^m\mu $ and ${\cal B}({\bf
K}^m)=\otimes_{i=1}^m{\cal B}({\bf K})$ is the minimal $\sigma $
algebra generated from the family of all subsets of the form
$A_1\times ... \times A_m$ with $A_i\in {\cal B}({\bf K})$ for each
$i=1,...,m$, where $U_{z,R}=B({\bf K}^m,z,R)$,
$x=z+e_1t_1+...+e_mt_m$, $t_i=t_{1,i}-t_{0,i}$, $\mbox{ }_iz=\mbox{
}_iu(t_{0,i})$, $z=(\mbox{ }_1z,...,\mbox{ }_mz)$. Thus $(1)$ is
satisfied for $\mu ^m$ almost all $z\in A^m$ and the applying of
Theorem 24.2 gives assertion 2.

\section{Appendix}
\par In \cite{sch1} it was proved the
non-archimedean variant of the Taylor theorem for functions of one
variable and for several variables in Theorem A.5 \cite{lujms04}.
\par {\bf A.1. Theorem.} {\it Let $f\in C^{n+1}(U,Y)$, where $n\in \bf
N$, $X$ and $Y$ be topological vector spaces over $\bf K$ and $U$ be
either clopen in $X$ or $X$ be locally convex. Then for each $x$ and
$y\in U$ the formula \\
$f(x)=f(y)+\sum_{j=1}^{n+1} {\bar {\Phi
}}^jf(y;x-y,...,x-y;0,...,0) +R_{n+1}(f;x,y).(x-y)^{\otimes (n+1)}$ \\
holds, where $R_{n+1}(f;x,y)=R_{n+1}(x,y): U^2\to L_{n+1}(X^{\otimes
(n+1)},Y)$ with \\ $\lim_{x\to y}R_{n+1}(x,y)=0$, $L_n(X^{\otimes
n},Y)$ denotes the space of $n$ polylinear continuous operators from
$X^{\otimes n}$ into $Y$.}
\par {\bf Proof.} If $j\le n$, then ${\bar {\Phi
}}^jf(z;v_1,...,v_j;0,...,0)$ is the $j$ polylinear operator by
vectors $v_1,...,v_j$ as follows from application of Lemma I.2. For
$n=0$ take $R_1(x,y).(x-y) := f(x)-f(y)-{\bar {\Phi
}}^1f(y;x-y;0)={\bar {\Phi }}^1f(y;x-y;1)-{\bar {\Phi
}}^1f(y;x-y;0)$. For $n=1$ from the definition of ${\bar {\Phi
}}^2f$ we have $f(x)-f(y)={\bar {\Phi }}^1f(y;x-y;1)$ and
$R_2(x,y).(x-y)^{\otimes 2} := {\bar {\Phi
}}^2f(y;x-y,x-y;0,1)-{\bar {\Phi }}^1f(y;x-y,x-y;0,0)$. Let the
statement be true for $n-1$, then from ${\bar {\Phi
}}^n(y+t_{n+1}(x-y);x-y,...,x-y;t_1,...,t_n)= {\bar {\Phi
}}^n(y;x-y,...,x-y;t_1,...,t_n)+ ({\bar {\Phi }}^1({\bar {\Phi
}}^nf(y;x-y,...,x-y;t_1,...,t_n))(y;x-y;t_{n+1}))t_{n+1}$ and the
continuity of ${\bar {\Phi }}^{n+1}f$ it follows that
$R_{n+1}(x,y).(x-y)^{\otimes (n+1)}={\bar {\Phi }}^nf(y;x-y,...,x-y;
0,...,0,1)-{\bar {\Phi }}^nf(y;x-y,...,x-y; 0,...,0,0)$, hence
$R_{n+1}(f;x,y)=R_{n+1}(x,y): U^2\to L_{n+1}(X^{\otimes (n+1)},Y)$.
Since $f\in C^{n+1}(U,Y)$, then $\lim_{x\to y}R_{n+1}(x,y)=0$. In
general the entire correction term $R_{n+1}(f;x,y).(x-y)^{\otimes
(n+1)}$ need not be polylinear, because $R_{n+1}(f;x,y)$ may be
nonlinear by $x-y$, but it is useful to write it in such form.
\par Here it is not used $D^i$, but partial difference quotients
${\bar {\Phi }}^i$ are used instead, then multipliers $1/i!$ does
not appear and the Taylor formula is true for $char ({\bf K})=p>0$
with decomposition up to terms of order $n+1\ge p$ as well.
Considering given $x, y\in U$ we can associate with it an
intersection of $\{ tx+(1-t)y: t\in {\bf K} \} $. If $U$ is clopen,
then $f$ has a $C^{n+1}$ extension on $X$, so we can consider a $\bf
K$ convex clopen subset $U_1$ such that $U\subset U_1\subset X$
instead of an initial one and denote it also by $U$. Thus, under
suppositions of this theorem on $U$ and $X$ we can consider the case
$\{ tx+(1-t)y: t\in B({\bf K},0,1) \} \subset U$ for each $x, y\in
U$ without loss of generality. Therefore,
$R_{n+1}(y+vt,y).(tv)^{\otimes
(n+1)}=t^{n+1}R_{n+1}(y+vt,y).v^{\otimes (n+1)} =({\bar {\Phi
}}^1({\bar {\Phi }}^n f(y;vt,...,vt;0,...,0))(y;v;t))t=t^{n+1}({\bar
{\Phi }}^1({\bar {\Phi }}^n f(y;v,...,v;0,...,0))(y;v;t))$, where
$x-y=vt$. Thus, the consideration can be reduced to $x$ and $y$
along lines containing $x$ and $y$. This gives \par
$R_{n+1}(y+vt,y).v^{\otimes (n+1)}=({\bar {\Phi }}^1({\bar {\Phi
}}^n f(y;v,...,v;0,...,0))(y;v;t))$,\\
where ${\bar {\Phi }}^n f(y;v_1,...,v_n;0,...,0)\in L_n(X^{\otimes
n},Y)$ and \par ${\bar {\Phi }}^1({\bar {\Phi }}^n
f(y;v_1,...,v_n;t_1,...,t_n))(y;v_{n+1};t_{n+1}))={\bar {\Phi
}}^{n+1} f(y;v_1,...,v_{n+1};t_1,...,t_{n+1})$\\ is continuous by
$(y;v_1,...v_{n+1};t_1,....,t_n,t_{n+1})\in U^{(n+1)}$ with
\par $\lim_{v_{n+1}\to 0}{\bar {\Phi }}^1({\bar {\Phi }}^n
f(y;v_1,...,v_n;0,...,0))(y;v_{n+1};t))=0$ \\ for each $t$ such that
$(y;v_1,...,v_{n+1};0,...,0,t)\in U^{(n+1)}$, since \\ ${\bar {\Phi
}}^{n+1}f(y;v_1,...,v_{n+1};0,...,0)\in L_{n+1}(X^{\otimes
(n+1)},Y)$ so that \\ ${\bar {\Phi
}}^{n+1}f(y;v_1,...,v_n,0;0,...,0)=0$, while ${\bar {\Phi
}}^{n+1}f(y;v_1,...v_n,0;t_1,...,t_n,t_{n+1})=0$ for $t_{n+1}\ne 0$
due to the definition of ${\bar {\Phi }}^{n+1}f$ and ${\bar {\Phi
}}^1f$. This proves the desired limit property of the residue
$R_{n+1}$ due to the continuity of ${\bar {\Phi }}^jf$ for each
$j\le n+1$, since the addition of vectors and multiplication of
vectors on scalars are continuous operations in $X$.

\end{document}